\documentclass[12pt,a4]{elsarticle}
\usepackage[
    pdftitle={IGA-HMM},
    pdfsubject={paper},
    pdfcreator={pdflatex},
    pdfstartview=FitV,
    pdfkeywords={iga, hmm, multiscale,   homogenization},
    bookmarksopen,
    pdfnewwindow=true,
    colorlinks=false,
    linkcolor=blue,
    citecolor=red,
    filecolor=magenta,
    urlcolor=cyan
    ]{hyperref}
\usepackage{hyperref}
\usepackage{fancyhdr}
\usepackage{tabularx}
\usepackage{graphicx}
\usepackage{subfigure}
\usepackage{epsfig}
\usepackage{booktaps}
\usepackage{float}
\restylefloat{figure}
\usepackage{amssymb,amsthm,makeidx,enumerate,amstext, amsmath,latexsym,amsbsy,euscript}
\usepackage{epstopdf}
\usepackage{natbib}
\usepackage{color}
\usepackage{threeparttable}

\makeatletter
\def\hlinewd#1{%
\noalign{\ifnum0=`}\fi\hrule \@height #1 %
\futurelet\reserved@a\@xhline}
\makeatother

\newcommand{\CC}{\mathbf{C}}

\newcommand{\PP}{\mathbf{P}}
\newcommand{\sS}{\mathbf{S}}

\newcommand{\xx}{\mathbf{x}}

\newcommand{\Eref}[1]{Eq.~(\ref{#1})}

\newcommand{\Fref}[1]{Fig.~\ref{#1}}

\newcommand{\cref}[1]{Chapter~\ref{chapt:#1}}

\newcommand{\cmatrixb}{\left\{ \begin{matrix}}
\newcommand{\cmatrixe}{\end{matrix} \right\}}

\newcommand{\vm}[1]{\mathbf{#1}}

\newcommand{\bc}{\begin{center}}
\newcommand{\ec}{\end{center}}
\newcommand{\bitem}{\begin{itemize}}
\newcommand{\eitem}{\end{itemize}}

\newcommand{\beq}{\begin{equation}}
\newcommand{\eeq}{\end{equation}}
\newcommand{\beqa}{\begin{eqnarray}}
\newcommand{\eeqa}{\end{eqnarray}}

\newcommand{\bv}{\begin{verbatim}}
\newcommand{\V}{\verb}                  % EX: \V=-d{#@~}= Expr must

\newcommand{\testpix}[1]{\fbox{\begin{minipage}[c]{\textwidth}
                      #1 \end{minipage} }}

\newcommand{\putpstex}[1]{\includegraphics{#1.pstex_t}}

\newcolumntype{C}{>{\centering\arraybackslash}X}
\begin{document}
\newtheorem{theorem}{Theorem}[section]
\newtheorem{lemma}[theorem]{Lemma}
\newtheorem{proposition}[theorem]{Proposition}
\newtheorem{corollary}[theorem]{Corollary}

\newtheorem{remark}[theorem]{Remark}

\newtheorem{definition}{Definition}[section] 
\newtheorem{example}{Example}[section]




\begin{frontmatter}
\title{An isogeometric analysis for elliptic homogenization problems}

\author[vn1,vn2]{H. Nguyen-Xuan\corref{nxh}}
\author[vn2]{T. Hoang}
\author[vn2]{V.P. Nguyen}

\cortext[nxh]{Corresponding author. Email address: nxhung@hcmus.edu.vn (H. Nguyen-Xuan).}

\address[vn1]{Department of Mechanics, Faculty of Mathematics and Computer Science, University of Science VNU HCMC, Ho Chi Minh City, Vietnam.}
\address[vn2]{Division of Computational Mechanics, Ton Duc Thang University, Ho Chi Minh City, Vietnam}

\begin{abstract}

A novel and efficient approach which is based on the framework of isogeometric analysis for elliptic homogenization problems is proposed. These problems possess highly oscillating coefficients leading to  extremely high computational expenses while using traditional finite element methods. The isogeometric analysis heterogeneous multiscale method (IGA-HMM) investigated in this paper is regarded as an alternative approach to the standard Finite Element Heterogeneous Multiscale Method (FE-HMM) which is currently an effective framework to solve these problems. The method utilizes non-uniform rational B-splines (NURBS) in both macro and micro levels  instead of standard Lagrange basis. Beside the ability to describe exactly the geometry, it tremendously facilitates high-order macroscopic/microscopic discretizations thanks to the flexibility of refinement and degree elevation with an arbitrary continuity level provided by NURBS basis functions. A priori error estimates of the discretization error coming from macro and micro meshes
and optimal micro refinement strategies for macro/micro NURBS basis functions of arbitrary orders are derived. Numerical results show the excellent performance of the proposed method.
\end{abstract}

\end{frontmatter}

\section{Introduction}

Homogenization is a branch of mathematics and engineering which studies partial differential equations (PDEs) with rapidly oscillating coefficients.  These kinds of equations describe various processes in inhomogeneous materials with rapidly oscillating micro structures, such as composite and perforated materials, and thus play an important role in physics, engineering and modern technologies. The aim of homogenization is to ``average out" the heterogeneities at micro-scale and describe the effective properties at macro-scale of such phenomena. In other words, we want to know the behaviors of the systems as homogeneous ones. In analytical approaches, such as in \citep{bensoussan1978,jikov_sm_1993}, homogenized equations are derived. However, the coefficients of these equations are only computed explicitly in some special cases, such as when the medium follows some periodic assumptions, and not explicitly available in general. Furthermore, full computations with complex scale interactions of the heterogeneous system are very ineffective due to high computational cost. Thus, to solve these problems, 
advanced computational technologies have been developed.

Literature reviews on various multiscale approaches can be found in \cite{weinan_principles_2011,multiscaletrends,PhuFE2Review}. In this paper, we focus on the heterogeneous multiscale method (HMM), which was proposed in \cite{WeinanEngquist2003}. Reviews on HMM are presented in \cite{WeinanHMMReview2007,AbdulleGakuto2009,AbdulleUnified2011,weinan_principles_2011}. This method provides a general framework which allows ones to develop various approaches to homogenization problems. The simplest one is the finite element heterogeneous multiscale method (FE-HMM), which uses standard finite elements such as simplicial or quadrilateral ones in both macroscopic and microscopic level. Solving the so-called micro problems (with a suitable set up) in sampling domains around traditional Gauss integration points at macro level allows one to approximate the missing effective information for the macro solver.

The standard finite element method (FEM), albeit very popular in various fields of sciences and engineering, still has some shortcomings which affect the efficiency of the FE-HMM. Firstly, the discretized geometry through mesh generation is required. This process often results in geometrical errors even with the higher-order FEM. Also, the communication between the geometry model and the mesh program during the analysis process is always needed and this constitutes the large part of the overall process \cite{HughesCottrellBazilevs:2005}, especially for industrial problems. Secondly, lower-order formulations, such as FE-HMM based on the four-node quadrilateral element (Q4), often require extremely fine meshes to produce approximate solutions with a desired accuracy for complicated problems. This prevents multiscale analyses from being run on personal computers.Thirdly, high-order formulations still put some restrictions on the element topologies (for example, the connection of different type of corner, center, or internal nodes) and only possess $C^0$ continuity. These disadvantages lead to an increase in the number of micro coupling problems and thus increase the overall computational cost. Hence, there is a need to consider alternative methods to tackle these issues.

Among advanced numerical methods, the so-called isogeometric analysis (IGA), where NURBS are used as basis functions, emerges as the most potential candidate. The isogeometric analysis was first proposed by Hughes and co-workers \cite{HughesCottrellBazilevs:2005} and now has attracted the attention of academic as well as  engineering community all over the world. The IGA provides a framework in which the gap between Computer Aided Design (CAD) and Finite Element Analysis (FEA) may be reduced. This is achieved in IGA by employing the same basis functions to describe both the geometry of the domain of interest and the field variables. While the standard FEM uses basis functions which are based on Lagrange polynomials, isogeometric approach utilizes more general basis functions such as B-splines and NURBS that are commonly used in CAD geometry. The exact geometry is therefore maintained at the coarsest level of discretization and re-meshing is seamlessly performed on this coarsest level without any further communication with CAD geometry. Furthermore, B-splines (or NURBS) provide a flexible way to make refinement, de-refinement, and degree elevation \cite{HughesCottrellBazilevs:2005}. They allow us to achieve easily the smoothness of arbitrary continuity in comparison with the $C^0$ continuity provided by the traditional FEM. For a reference on IGA, we recommend the excellent book \cite{CottrellHughesBazilevs:2009} and we refer to the NURBS book \cite{piegl_nurbs_1997} for a geometric description.


In this work, we introduce a new approach: a so-called isogeometric analysis heterogeneous multiscale method (IGA-HMM) which utilizes NURBS as basis functions for both exact geometric representation and analysis. The NURBS are used as basis functions for both macro and micro element spaces, where the  FE-HMM only employs standard FEM basis. This tremendously facilitate high-order macroscopic/microscopic discretizations by a flexibility of refinements and degree elevations with an arbitrary continuity of  basis functions. As will be demonstrated later in the  numerical examples section, 
elliptic homogenization problems can be solved, by using the proposed IGA-HMM, effectively on a personal computer which is not the case for FE-HMM with (bi)linear basis functions if high accuracies are needed. 
The FE-HMM with bi(linear) elements often requires very fine macro meshes (thus a high number of micro problems) that are only supported by powerful computers. It is obvious that one can use high order Lagrange elements in the framework of FE-HMM to achieve the same goal. However, high order Lagrange elements cannot be constructed as straightforwardly as NURBS and more importantly they are only $C^0$ continuity whereas NURBS elements
are $C^{p-1}$ ($p$ is the NURBS order) by definition. Therefore, IGA-HMM is able to solve high order homogenization problems such as plate/shell homogenization problems. One can consider our IGA-HMM as an efficient high order NURBS-based FE-HMM. We refer to \cite{li_efficient_2012} for a related work on high order FE-HMM which, however, limits to cubic macro elements and quadratic micro elements. A priori error estimates of the discretization error coming from macro and micro meshes are  provided. Optimal micro refinement strategies for macro and micro NURBS basis functions of arbitrary orders are presented and thus an optimal value for the micro NURBS order is obtained as a function of the macro NURBS basis order. 
We have also observed that $C^0$, not $C^{p-1}$, high order NURBS should be used at the microscale.
It should be mentioned that the proposed IGA-HMM is very similar to the standard FE-HMM, this means that the implementation is very simple, any existing coding framework of FE-HMM e.g. the one given in \citep{AbdulleShort2009} can be readily re-used.

The paper is arranged as follows: a brief introduction to the B-splines, NURBS and IGA is given in Section 2. Section 3 describes the isogeometric analysis heterogeneous multiscale method. Numerical examples are provided in Section 4 and Section 5 closes the paper with some concluding remarks.

\section{NURBS-based isogeometric analysis fundamentals}

\subsection{Knot vectors and basis functions}
Univariate B-spline basis functions are constructed from
\textit{knot vectors}. Given two positive integers $p$ and $n$ and
let $\Xi= \left[\xi_1,\xi_2,...,\xi_{n+p+1}\right]$ be a
non-decreasing sequence of parameter values, ${\xi_i} \le
{\xi_{i + 1}},i = 1,...,n+p$.
The $\xi_i$ are called knots, $\Xi$ is the set of coordinates in
the parametric space, $p$ is the polynomial degree, and n is the
number of basis functions. If all knots are equally spaced the
knot vector is called uniform.
Otherwise, they are called a non-uniform knot vector. When the
first and the last knots are repeated $p +1$ times, the knot
vector is called open. An important property of open knot vectors is that the resulting basis functions are interpolatory at the ends of the parametric space.
A B-spline basis function is $C^\infty$
continuous inside knot spans and $C^{p-1}$ continuous at a
single knot. If an interior knot value repeats one more time, it
is then called a multiple knot. At a knot of multiplicity $k$
the continuity is $C^{p-k}$.
\noindent Given a knot vector, the B-spline basis functions
$N_{i,p}(\xi)$ are defined starting with the zeroth order
$(p=0)$ basis functions

\begin{equation}
{N_{i,0}}(\xi) = \left\{ {\begin{array}{*{20}{c}}
   {1 ,~~ \text{if} ~~{\xi_i} \le \xi \le {\xi_{i + 1}}},\\
   {0,~~\text{otherwise;~~~~~~~~}}  \\
\end{array}} \right.
 \label{eq:equation17}
 \end{equation}
and for $p\geq 1$, they are defined recursively as follows
\cite{piegl_nurbs_1997}
\begin{equation}
{N_{i,p}}\left( \xi \right) = \frac{{\xi - {\xi _i}}}{{{\xi _{i
+ p}} - {\xi _i}}}{N_{i,p - 1}}\left( \xi \right) + \frac{{{\xi
_{i + p + 1}} - \xi }}{{{\xi _{i + p + 1}} - {\xi _{i +
1}}}}{N_{i + 1,p - 1}}\left( \xi \right).
\end{equation}
For $p = 0$ and $p=1$, the NURBS basis functions of isogeometric
analysis are identical to those of standard piecewise constant
and linear finite elements, respectively. Nevertheless, for
$p\geq 2$, they are different
~\citep{HughesCottrellBazilevs:2005}. Therefore, the present
work will consider the basis functions with $p\geq 2$. 

\Fref{fig:NURBS_basis} illustrates a set of univariate 
quadratic, cubic and quartic B-spline basis functions for open
uniform knot vectors $\Xi =\{0,0,0,\frac{1}{2},1,1,1\}$, $\Xi
=\{0,0,0,0,\frac{1}{2},1,1,1,1\}$ and $\Xi
=\{0,0,0,0,0,\frac{1}{2},1,1,1,1,1\}$, respectively.

\begin{figure}[!ht]
\vspace{-.5cm}
\centering
\subfigure{\includegraphics[width=0.32\columnwidth]{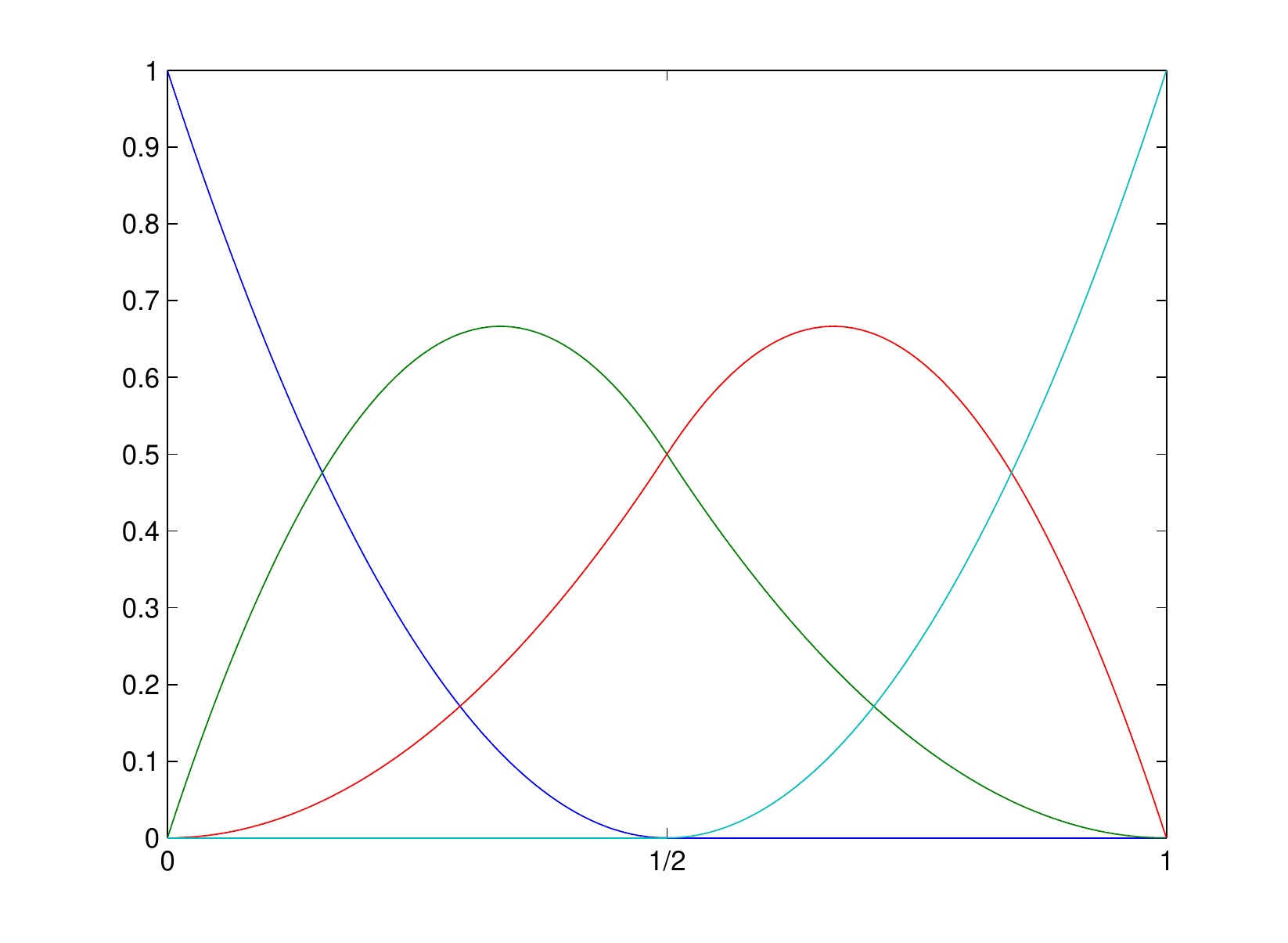}}
\subfigure{\includegraphics[width=0.32\columnwidth]{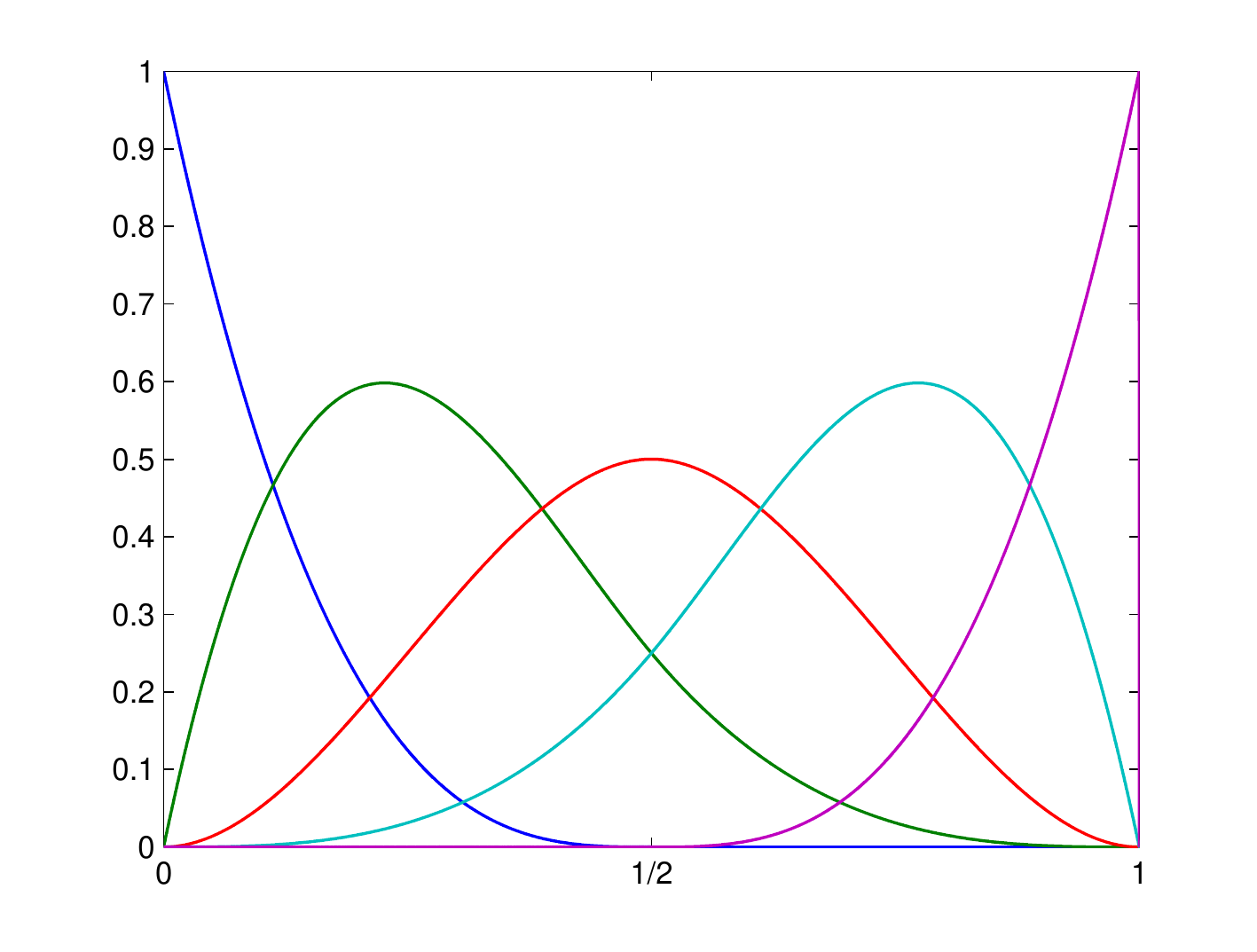}}
\subfigure{\includegraphics[width=0.32\columnwidth]{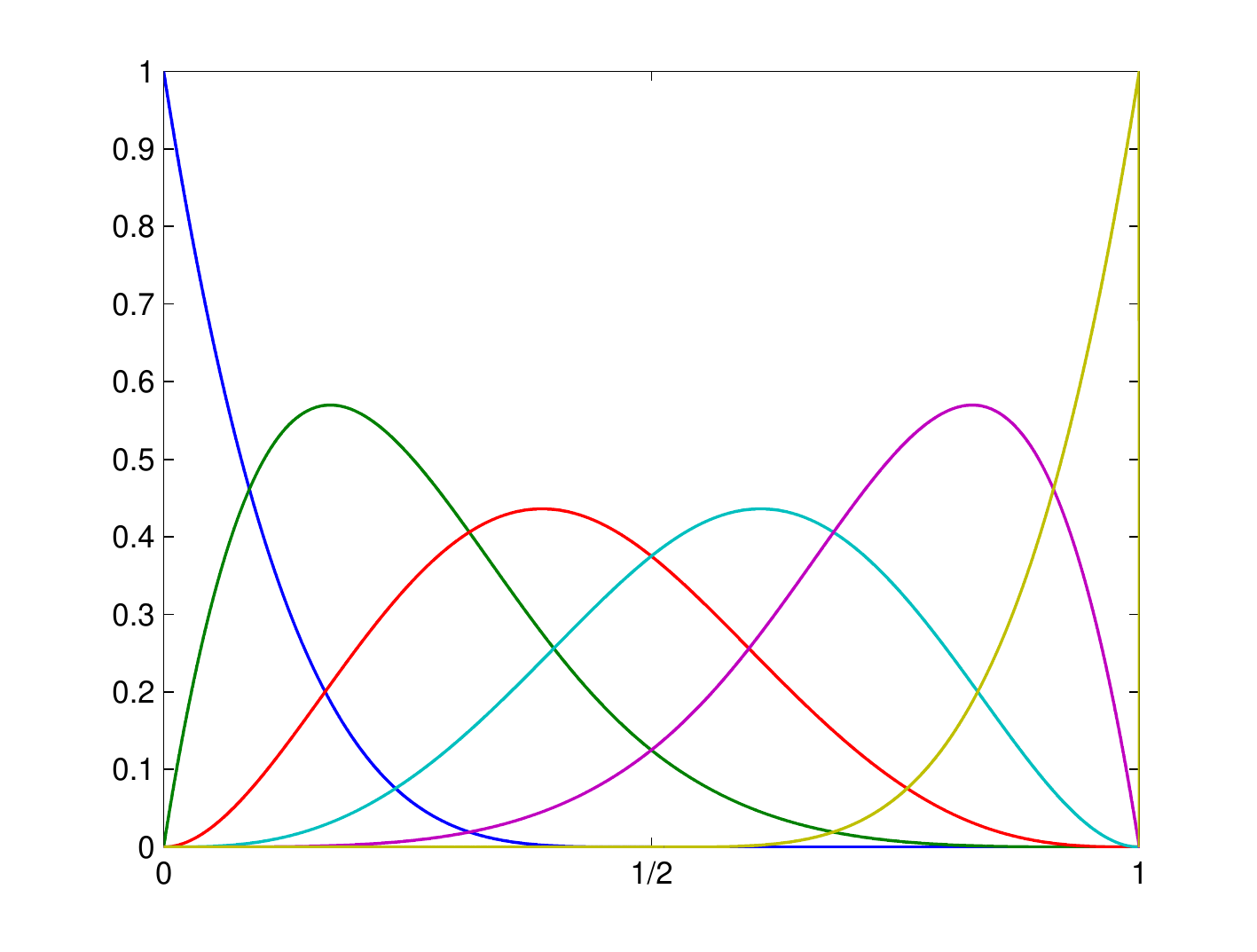}}
\vspace{-.5cm}
\caption[This figure illustrate cubic shape basis of B-spline]
{An illustration of quadratic, cubic and quartic B-spline basis functions.}
\label{fig:NURBS_basis}
\end{figure}

%

\subsection{NURBS curves and surfaces}
The B-spline curve is defined as
\begin{equation}
\CC\left( \xi \right) = \sum\limits_{i = 1}^n {{N_{i,p}}} \left(
\xi \right){\PP_i},
\end{equation}
where $\PP_i$ are the control points and $N_{i,p}\left( {\xi}
\right)$ is the $p^{th}$-degree B-spline basis function defined
on an open knot vector. 

The B-spline surfaces  are defined by the tensor product of
the basis functions in two parametric dimensions $\xi$ and
$\eta$ with two knot vectors $\Xi = \left\{ {{\xi _1},{\xi
_2}...,{\xi _{n + p + 1}}} \right\}$ and $\mathcal{H}= \left\{
{{\eta _1},{\eta _2}...,{\eta _{m + q + 1}}} \right\}$ as

\begin{equation}
\sS\left( {\xi ,\eta } \right) = \sum\limits_{i = 1}^n
{\sum\limits_{j = 1}^m {{N_{i,p}}\left( \xi
\right){M_{j,q}}\left( \eta \right){\PP_{i,j}}} },
\label{eq:equation1}
\end{equation}
where $\PP_{i,j}$ is the bidirectional control net,
$N_{i,p}\left( {\xi} \right)$ and $M_{j,q}\left( {\eta }
\right)$ are the B-spline basis functions defined on the knot
vectors over an $n\times m$ net of control points $\PP_{i,j}$.
Following the traditional FEM notation, \Eref{eq:equation1} is
rewritten as follows
\begin{equation}
\sS\left( {\xi ,\eta } \right) = \sum\limits_I^{n\times m} {{N^b_I}}
\left( {\xi ,\eta } \right){\PP_I},
\end{equation}
where ${N^b_I}\left( {\xi ,\eta } \right) = {N_{i,p}}\left( \xi
\right){M_{j,q}}\left( \eta \right)$ is the bivariate B-spline basis function associated with node I.
Fig. \ref{fig:2d_NURBS_basis} plots the bivariate quadratic and cubic B-spline basis functions.

\begin{figure}[!ht]
\begin{center}
\subfigure{\includegraphics[angle=0,width=0.38\columnwidth]{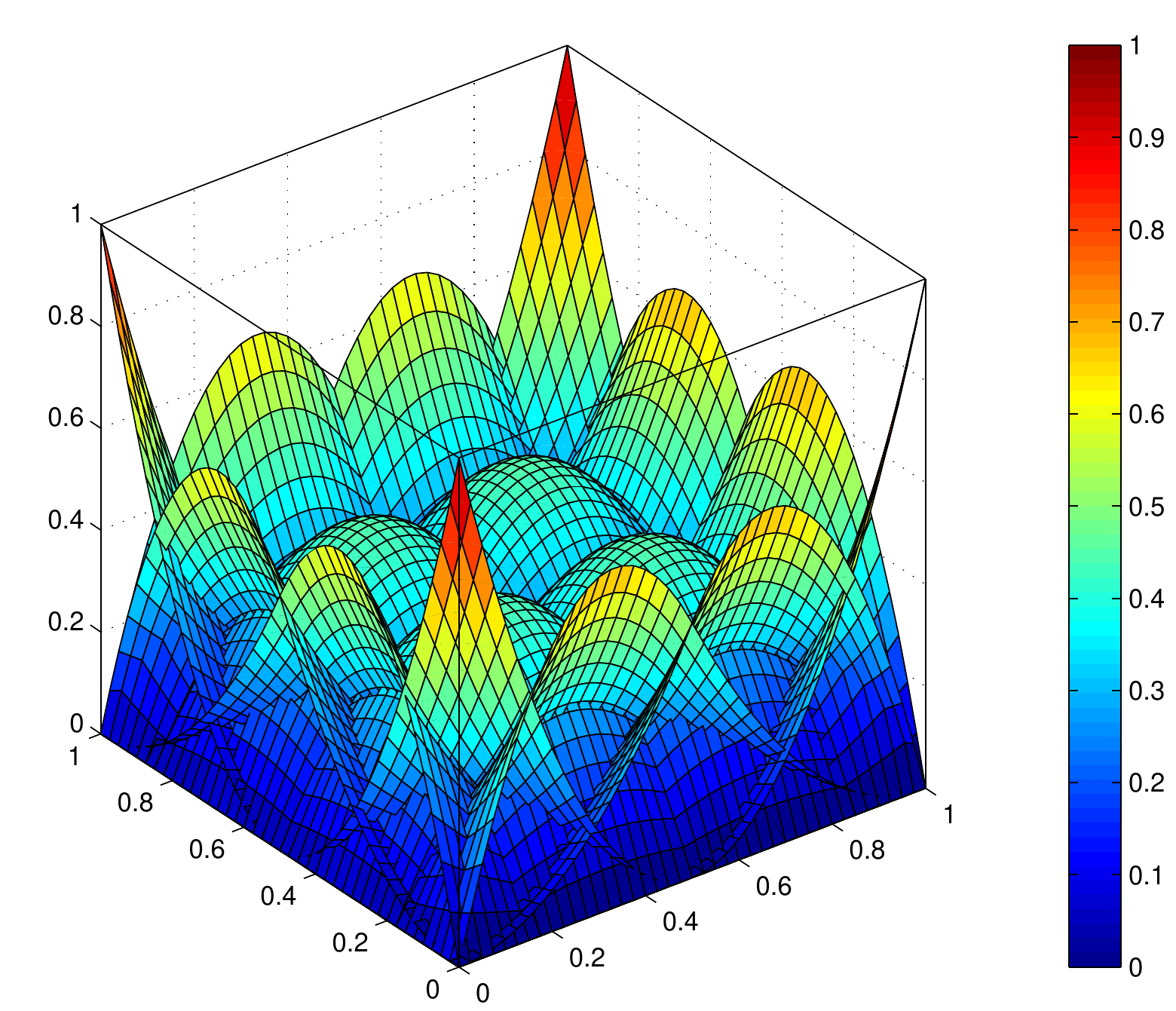}}
\subfigure{\includegraphics[angle=0,width=0.38\columnwidth]{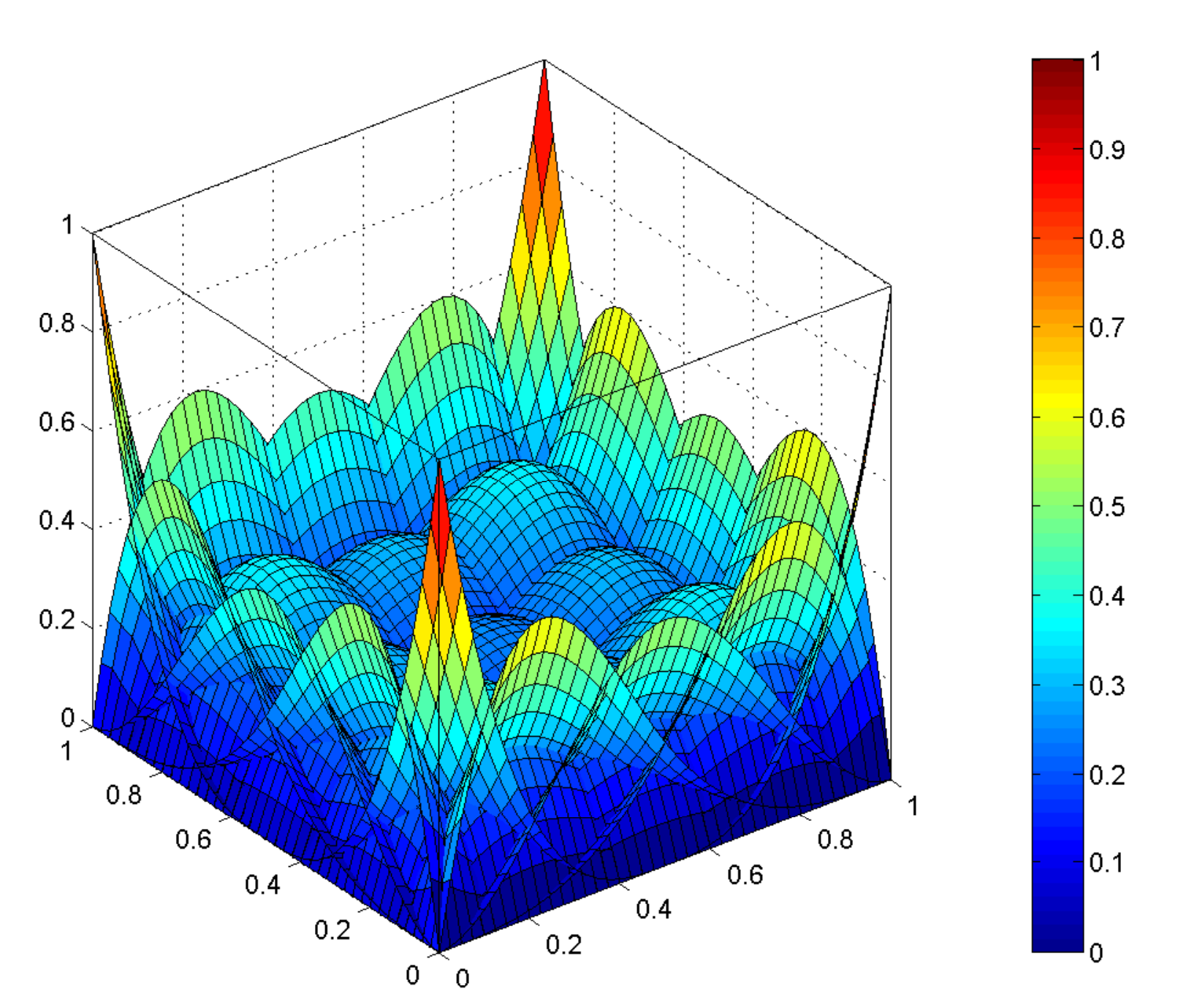}}
\end{center}
\vspace{-.5cm}
\caption {Bivariate quadratic and cubic B-spline basis functions with open uniform 
knot vectors $\Xi =\mathcal{H}=\{0,0,0,\frac{1}{2},1,1,1\}$, $\Xi=\mathcal{H}
=\{0,0,0,0,\frac{1}{2},1,1,1,1\}$.}
\label{fig:2d_NURBS_basis}
\end{figure}

B-splines are convenient for free-form modelling, but they lack the ability to exactly represent some simple shapes like circles and ellipsoids.
Non-uniform rational B-splines (NURBS) extend B-splines since they allow exact representation of conic sections.
NURBS are obtained by
augmenting every point in the control mesh $\PP_I$ with the weight $w_I$. The weighting function is
constructed as follows
\begin{equation}
w\left( {\xi ,\eta } \right) = \sum\limits_{I = 1}^{n\times m}
{{N^b_I}} \left( {\xi ,\eta } \right){w_I}.
\end{equation}
\noindent The NURBS surfaces are then defined by

\begin{equation}
\sS\left( {\xi ,\eta } \right) = \frac{{\sum\limits_I^{n\times m}
{{N^b_I}} \left( {\xi ,\eta } \right){w_I}{\PP_I}}}{{w\left( {\xi
,\eta } \right)}} = \sum\limits_{I = 1}^{n\times m} {{{
N}_I}\left( {\xi ,\eta } \right)} {\PP_I},
\end{equation}
where ${{N}_I}\left( {\xi ,\eta } \right) =
\frac{{{N^b_I}\left( {\xi ,\eta } \right){w_I}}}{{w\left( {\xi
,\eta } \right)}}$ are the rational basis functions. An example of quadratic NURBS curve and surface is
demonstrated in \Fref{fig:NURBS_surface}. It is clear that, due to the partition of unity property of B-splines basis, if all the weights are units then NURBS curves/surfaces become B-splines curves/surfaces.

\begin{figure}[H]
\begin{center}
\subfigure{\includegraphics[angle=0,width=0.49\columnwidth]{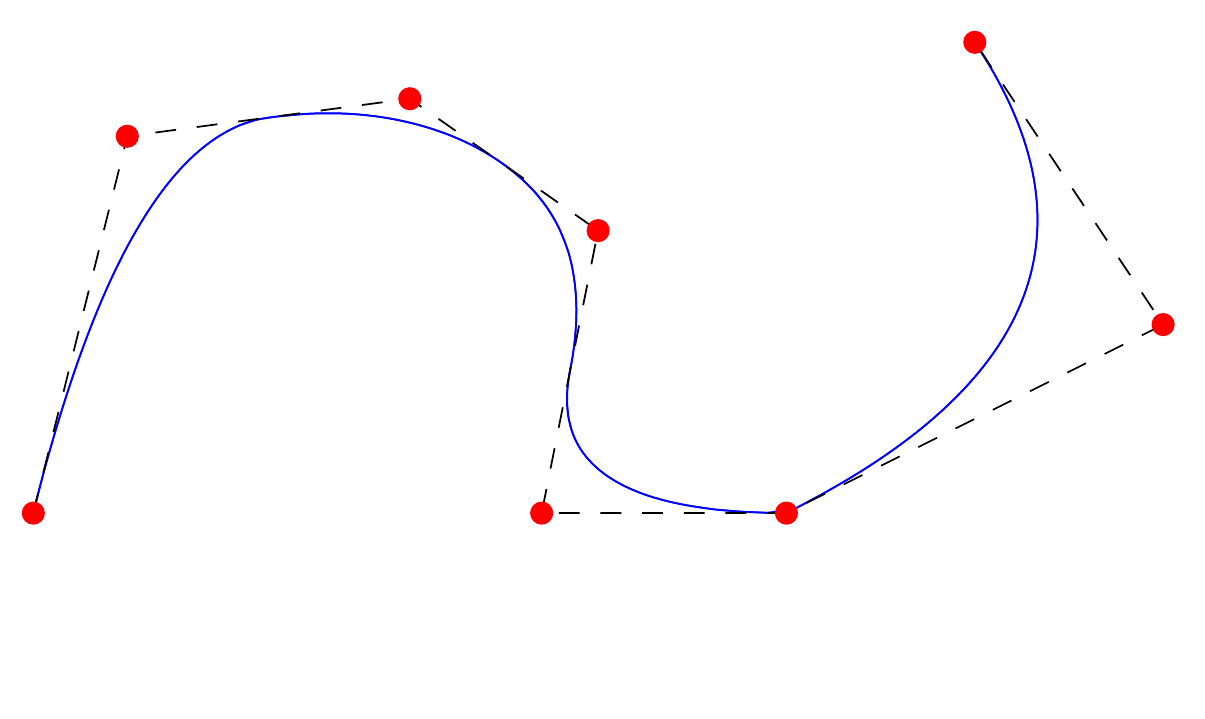}}
\subfigure{\includegraphics[angle=0,width=0.49\columnwidth]{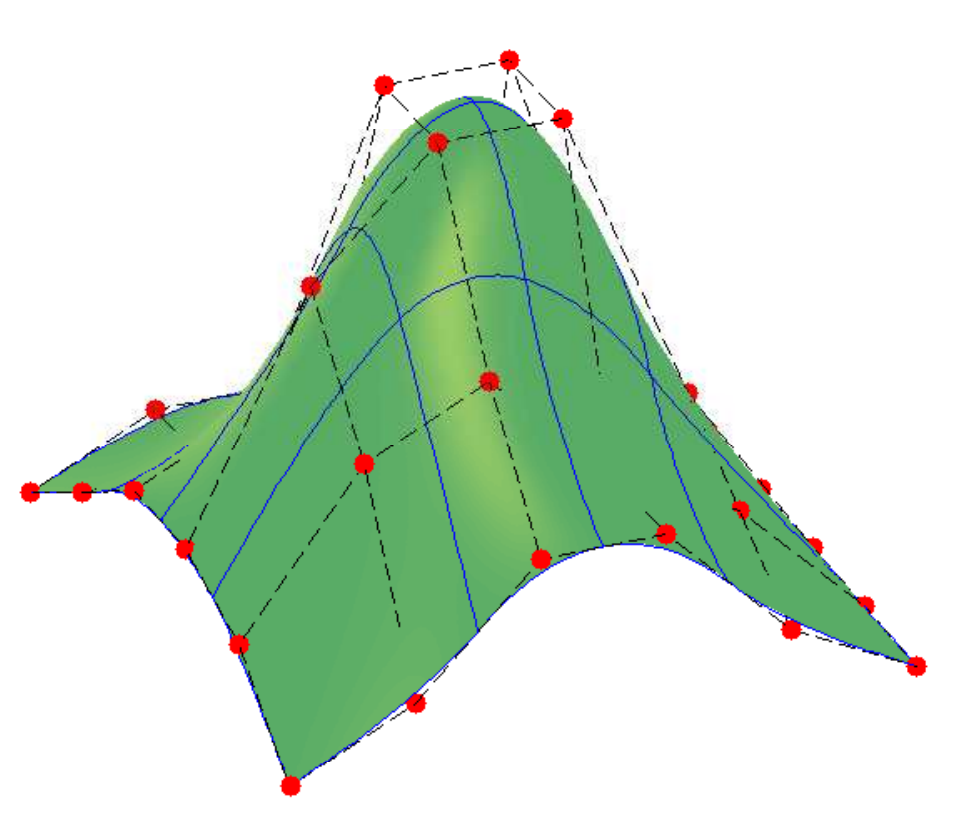}}
\end{center}
\caption{{\color{black}A quadratic NURBS curve (left) and a quadratic NURBS surface  with physical mesh and control grid (right)}}
 \label{fig:NURBS_surface}
\end{figure}

\subsection{Refinement}

In addition to the classical $h$-refinement (knot insertion in CAD notation) 
and $p$-refinement (degree elevation in CAD notation) supported by FEM, IGA has a new, unique, and more economical
version: $k$-refinement which is a process in which degree elevation is followed by knot insertion. The $k$-refinement allows us to
elevate the degree of NURBS objects while also obtaining higher
continuities across element boundaries, whereas the standard FEM
always stay at $C^0$. In comparison to $p$-version, $k$-version
require much less number of basis functions while maintaining
the same convergence rate. The readers are referred to
\cite{HughesCottrellBazilevs:2005, CottrellHughesBazilevs:2009} for more details.

\subsection{Isoparametric concept}

Isogeometric analysis also employs the isoparametric concept as in standard FEM--the same NURBS basis functions are used for both the description of the geometry and the unknown field:

\begin{equation}
\xx (\xi,\eta)= \sum\limits_I^{n\times m} {{N_I}\left( {\xi ,\eta } \right)} {\PP_I},~~~~
u^H(\xx (\xi,\eta))= \sum\limits_I^{n\times m} {{N_I}\left( {\xi ,\eta } \right)} {u_I},
\label{eq:isoparametric}
\end{equation}
where $n\times m$ represent the number of basis functions, $\xx^T=(x,y)$ is the physical coordinates  vector, ${N_I}\left( {\xi ,\eta }\right)$ is the NURBS basis function and $u_I$ is the degree of freedom (dof) of $u^H$ (which can be the displacements for solid mechanics problems, the temperatures for thermal problems etc.) at the control point $I$, respectively. The control points play the role of nodes in a FEM setting and the IGA elements are defined as the non-zero knot spans. 

\subsection{Properties of IGA}

The NURBS-based IGA and FEA share a lot of commons such as \cite{HughesCottrellBazilevs:2005}:  compact support, partition of unity, affine covariance, isoparametric concept, patch test satisfied. These properties allow implementations in the IGA to follow the framework of Galerkin method as in the FEA. On the other hand, there are important properties which IGA possesses but FEA does not. Some differences between NURBS-based IGA and FEA are presented in Tab. \ref{diffIGAFEA}. Note that one weak point of NURBS-based IGA is the lack of interpolation property which renders the application of Dirichlet boundary conditions  a bit more involved than traditional FEM. We refer to \cite{HughesCottrellBazilevs:2005} for details on this issue. One drawback of NURBS is the lack of local refinement, thus making adaptive NURBS-based IGA impossible. To this end, T-splines have been recently developed to overcome the limitations of NURBS while making use of existing NURBS algorithms
\cite{bazilevs_isogeometric_2010}. Another solution is the hierarchical refinement approach presented in \cite{vuong_hierarchical_2011}.
These two options can be straightforwardly incorporated into our IGA-HMM.

\begin{table}[H]
\caption{Differences between NURBS-based Isogeometric analysis and Finite element analysis \cite{HughesCottrellBazilevs:2005}}
\begin{center}
\begin{tabular}{l l l}
\hlinewd{1.2pt}
& IGA  & FEA\\ \hlinewd{0.8pt}
Geometry & exact  & approximation\\
Handling points & control points  & nodal points\\
Variables & control variables & nodal variables\\
Basis & NURBS  & Lagrange \\
Basis interpolation property & no & yes\\
Continuity & easily controlled & $C^0$, fixed \\
Refinement space & $hpk$ & $hp$\\
Positive basis  property& yes  & no\\
Convex hull property & yes & no \\
Variation dismishing property & yes  & no\\
\hlinewd{1.2pt}
\end{tabular}
\label{diffIGAFEA}
\end{center}
\end{table}

Fig. \ref{fig:CkC0} plots two meshes that consist of quadratic NURBS elements--one with $C^1$ continuity across element edges (Fig. \ref{fig:CkC0}a) and one with only $C^0$ continuity across element edges (Fig. \ref{fig:CkC0}b). The latter was seamlessly 
obtained from the former by using the so-called knot insertion operation. Also clear from the referred figure is the fact that the number of
dofs of the $C^1$ mesh is less than the number of dofs of the $C^0$ mesh. It is worthy noting that as 
a result of several decades of research, many efficient computer algorithms exist for the fast evaluation and refinement
of NURBS  \cite{piegl_nurbs_1997}. Therefore, NURBS-based FEM is very efficient.

\begin{figure}
        \centering
        \subfigure[$C^1$ quadratic NURBS mesh]{\includegraphics[width=0.49\textwidth]{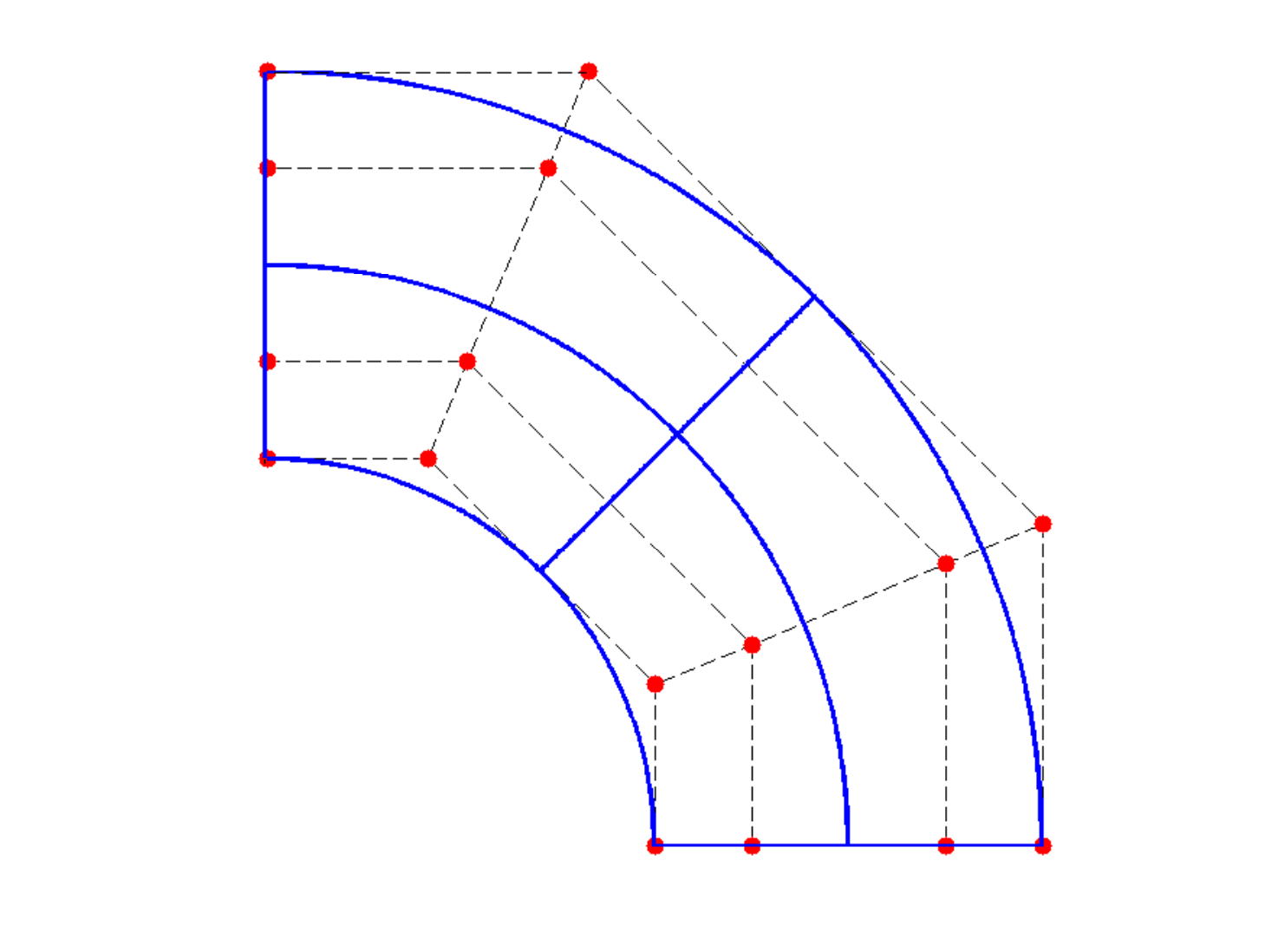}}
        \subfigure[$C^0$ quadratic NURBS mesh]{\includegraphics[width=0.49\textwidth]{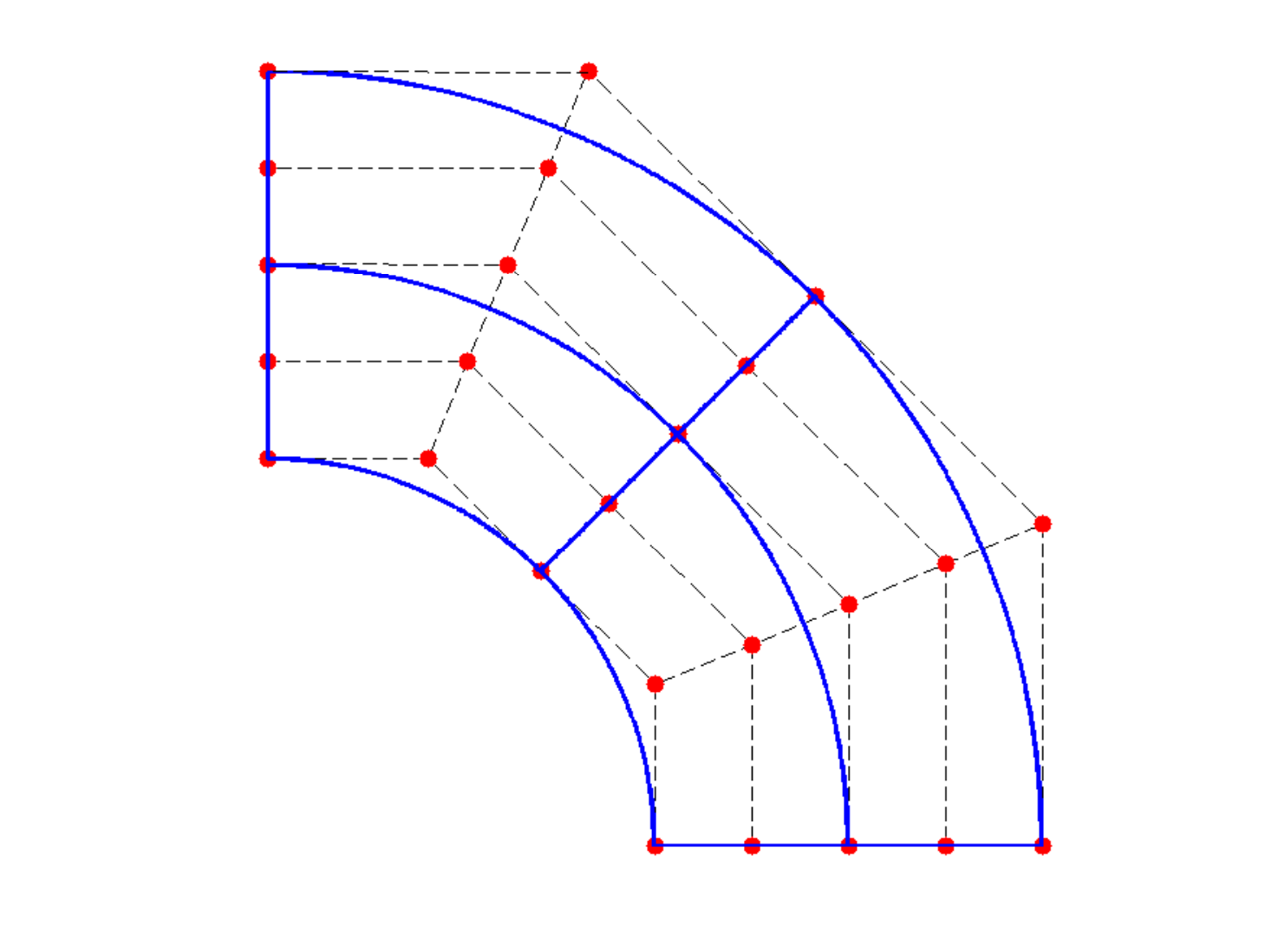}}
        \caption{Quadratic NURBS meshes ($2\times2$ elements): $C^1$ case (a) and $C^0$ case (b).
        Red dots denote the control points which generally do not locate inside the physical domain.}
        \label{fig:CkC0}
\end{figure}


\section{An isogeometric analysis heterogeneous multiscale method (IGA-HMM)}

\subsection{Model problems}

Let $\Omega $ be a domain in ${{\mathbb{R}}^{d}}(d=1,2,3)$ with
Lipschitz boundary $\partial \Omega$ on which we impose the
Dirichlet boundary conditions. Given $f\in {{L}^{2}}\left( \Omega
\right)$ as a source term, we consider the following classical
highly oscillating coefficient second-order elliptic problem

\begin{equation}\label{eq1}
\left\{ \begin{matrix}
-\nabla \cdot \left( {{a}^{\varepsilon }}(\mathbf{x})\nabla
{{u}^{\varepsilon }}(\mathbf{x}) \right) & = & f(\mathbf{x}) & \text{in }\Omega ,
\\
{{u}^{\varepsilon }}(\mathbf{x}) & = & 0, & \text{on }\partial \Omega,
\\
\end{matrix} \right .
\end{equation}
where

\begin{itemize}
\item conductivity tensor ${{a}^{\varepsilon }}\in {{\left(
{{L}^{\infty }}\left( \Omega \right) \right)}^{d\times d}}$ ,
and is uniformly elliptic and bounded,
\item $\varepsilon $ is a parameter which represents a fine
scale characterizing the multiscale nature of ${{a}^{\varepsilon
}}$.
\end{itemize}
and without loss of generality, we have omitted the Neumann boundary conditions.

Equation \eqref{eq1} is popular in many practical problems such as solid mechanics, thermal conduction, and electrostatics, etc, where $u$ can be the displacement field or electric field, etc, respectively, see Tab. \ref{significant}.

\begin{table}[H]
\caption {Significance of $u$ and $f$ in application \cite{felippa:note}}
\begin{center}
\begin{tabular}{llll}
\hlinewd{1.2pt}
Application Problem& $a^\varepsilon$ & $u^\varepsilon$  &  $f$ \\
\hlinewd{0.8pt}
Solid mechanics &Material tensor& Displacement & Mechanical force\\
Heat conduction &Thermal conductivity& Temperature  & Heat flux\\
Acoustic fluid &Acoustic conductivity&  Displacement potential  & Particle velocity\\
Potential flows &Flow conductivity&  Pressure  & Particle velocity\\
General flows  &Flow conductivity &Velocity  & Fluxes\\
Electrostatics  &Electrical conductivity& Electric potential  & Charge density\\
Magnetostatics  &Magnetical conductivity& Magnetic potential  & Magnetic intensity\\
\hlinewd{1.2pt}
\end{tabular}\label{significant}
\end{center}
\end{table}

In this paper, we confine ourselves to two dimensional heat elliptic problems, where $u^\varepsilon$  is the temperature field.  We emphasize that the applications for other problems are similar and extension to three dimensions is straightforward \citep{AbdulleShort2009}.

From homogenization theory, we have known that the solution
${{u}^{\varepsilon }}$ of (1) (weakly) converges to ${{u}^{0}}$,
which is the solution of a so-called \textit{homogenization
problem}
\begin{equation}
\label{eq2}
\left\{ {\begin{array}{*{20}{c}}
  { - \nabla \cdot\left( {{a^0}({\mathbf{x}})\nabla {u^0}({\mathbf{x}})} \right)}& = &{f({\mathbf{x}}),}&{{\mathbf{x}} \in \Omega  \subset {\mathbb{R}^d},} \\
  {{u^0}({\mathbf{x}})}& = &{0,}&{{\mathbf{x}} \in \partial \Omega ,}
\end{array}} \right.{\text{ }}
\end{equation}
where the so-called homogenized tensor ${{a}^{0}}$ is usually
not available nor explicitly analytically computed.

The FE-HMM or IGA-HMM aims at finding the upscale solution ${{u}^{0}}$ without computing $a^0$ explicitly.

\subsection{Drawbacks of the FE-HMM method}
One of the shortcomings of the FE-HMM is degree elevation. Due to the connection of different types of corner, center, or internal
nodes, the implementation for degree elevation of standard FEM basis
functions is, albeit possible, not an easy task. So far, to the best of the authors' knowledge,  in the FE-HMM numerical
literature, the highest order basis function has been used is cubic
\cite{li_efficient_2012}. Another drawback of FE-HMM is that 
geometry errors will appear when applied to curved boundary domains.
This is because in FEA meshing, polygons are used to approximately
represent curved boundaries. To reduce this geometry error, we have
to rely on sufficiently fine meshes, but this will lead to a high
number of micro problems. These issues can be seamlessly solved when we use a
new approach--the so-called isogeometric analysis heterogeneous multiscale
method (IGA-HMM). Moreover, the IGA-HMM is able to utilize $k$-refinement - higher order and higher continuity, which is a
unique advantage of IGA \cite{CottrellHughesBazilevs:2009}.


\subsection{The isogeometric analysis heterogeneous multiscale method (IGA-HMM)}

The crucial idea of the IGA-HMM is
to replace the standard FEM basis functions used in FE-HMM by the
NURBS basis ones, in both macro and micro scales. Before
describing the IGA-HMM method, we introduce the so-called macro
and micro finite ``patch'' spaces in which we will work.

\renewcommand{\thefootnote}{\fnsymbol{footnote}}
\setcounter{footnote}{1}
Let $\Omega$ be an open subset of $\mathbb{R}^d, d=1, 2, 3,$ and denote $Y =  ~(0,1)^d $ the unit cube in  $\mathbb{R}^d$. We list here some familiar notations:\\
\begin{equation*}
\begin{split}
 {L^2}\left( \Omega  \right) = &~\left\{ {v:\Omega  \to \mathbb{R}:\int_\Omega  {{{\left| v \right|}^2}d\mathbf{x}}  < \infty } \right\},  ~{\left\| v \right\|_{{L^2}\left( \Omega  \right)}} = ~\left(\int_\Omega  {{{\left| v \right|}^2}d\mathbf{x}}\right)^{1/2} ;\hfill \\
  {H^1}\left( \Omega  \right) =&~\left\{ {v \in {L^2}\left( \Omega  \right):\nabla v \in {\left( {{L^2}(\Omega )} \right)^d} } \right\}, ~{\left\| v \right\|_{{H^1}\left( \Omega  \right)}} = ~{\text{ }}{\left( {\int_\Omega  {({{\left| v \right|}^2} + {{\left| {\nabla v} \right|}^2)}d\mathbf{x}} } \right)^{1/2}};  \hfill \\
   H_0^1\left( \Omega  \right) = &~\left\{ {v \in {H^1}\left( \Omega  \right):v = 0\,{\text{on}}\,\partial \Omega } \right\}; \hfill \\
  H_{per}^1(Y) = &~\left\{ {{v \in {H^1 }(Y):v \text{ is } Y-periodic \footnotemark}}\right\}.\hfill  \\
 \end{split}
\end{equation*}
\footnotetext{For $Y = (0,1)^d $, a function $f:{\mathbb{R}^d} \to \mathbb{R}$ is said to be $Y-periodic$ if and only if $f(\mathbf{x}+\textbf{e}_j)=f(\mathbf{x}) \text{ for all } \mathbf{x} \in {\mathbb{R}^d}, j\in \{1,..,d\}$, where $\{\textbf{e}_j\}_{j=1}^{d}$ is the canonical basis of ${\mathbb{R}^d}$. Thus we may view a function $v \in H_{per}^1(Y) $ as a function from $Y$ to ${\mathbb{R}}$ which belongs to $H^1(Y)$ and has the same trace on opposite faces of $Y$. }

\begin{remark}
In this work, we assume NURBS with equal orders in two directions ($p=q$). At the macroscale, NURBS of order $p$ and at the microscale,
NURBS of order $q$ are utilized. This facilitates the derivation of error estimations presented later.
\end{remark}

\begin{figure}[H]
\begin{center}
\subfigure{\includegraphics[trim= 2cm 1cm 2cm 1cm,
angle=0,width=0.75\columnwidth]{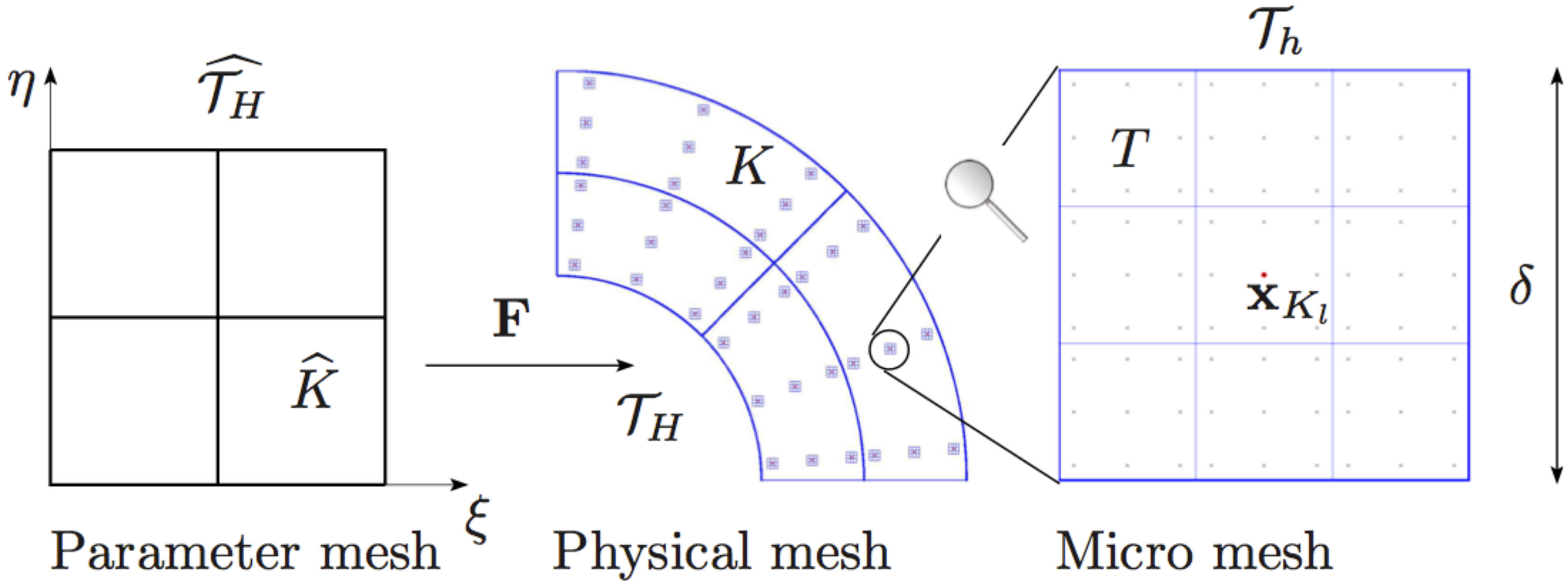}}
\end{center}
\caption{IGA-HMM: from parameter mesh to physical mesh at the macro scale. For every Gauss point $\vm{x}_{K_l}$ of the macro element $K$, 
a micro domain $K_{\delta_l}$ is considered. The triangulation of $K_{\delta_l}$ is ${{\mathcal{T}}_{h}}$. Solution of the micro problem
using the micro mesh ${{\mathcal{T}}_{h}}$ yields the missing information at the macro Gauss point $\vm{x}_{K_l}$.}
\label{fig:hmm}
\end{figure}

\subsubsection{Macro finite patch space}

We denote by $\Xi =\left\{ {{\xi }_{1}},...,{{\xi }_{n+p+1}}
\right\},\mathcal{H}=\left\{ {{\eta }_{1}},...,{{\eta }_{m+p+1}}
\right\}$ two knot vectors which define the index space, and
$\textbf{F}$ is a global geometry map function which
parameterizes the physical space $\Omega $
\[\textbf{F}:\widehat{\Omega }\to \Omega, \]
$\textbf{F}$ is defined by the first equation of Eq. \ref{eq:isoparametric}.
The knot vectors $\Xi $ and $\mathcal{H}$ define, on the parameter
space $\widehat{\Omega }$, a macro parameter mesh
$\widehat{{{\mathcal{T}}_{H}}}$. Through the geometry map $\textbf{F}$, the parameter mesh
$\widehat{{{\mathcal{T}}_{H}}}$ become the physical mesh
${{\mathcal{T}}_{H}}$ on the domain $\Omega $. We denote by
${{H}_{K}}$ the diameter of elements $K\in {{\mathcal{T}}_{H}}$,
and set $H:={{\max }_{K\in {{\mathcal{T}}_{H}}}}{{H}_{K}}$. We refer to Fig. \ref{fig:hmm} for a graphics
illustration.

The macro finite patch space is then defined as

\begin{equation}
S_{0}^{p}\left( \Omega ,{{\mathcal{T}}_{H}}
\right)=H_{0}^{1}\left( \Omega \right)\cap \text{span}{{\left\{
R_{ij}^{{}}\circ {{\textbf{F}}^{-1}} \right\}}_{i,j}},
\end{equation}
where $\{{R}_{ij}\}$ are the two dimensional NURBS basis functions
of order $p$ defined on the parameter space $\widehat{\Omega }$.

On each element $\widehat{K}\in \widehat{{{\mathcal{T}}_{H}}}$,
we consider ${{\left\{ {{\omega
}_{{{K}_{l}}}},\widehat{{{\mathbf{x}}_{{{K}_{l}}}}}
\right\}}_{l=1,..,(p+1)^2}}$ the Gauss quadrature weights and
integration points\footnote{In 2D, these Gauss points are defined as follows. First, they are defined in the usual parent
domain $[-1,-1]\times[1,1]$. Next, they are transformed to $\widehat{K}$. Note that this step does not exist in a standard FEM. Refer to \cite{CottrellHughesBazilevs:2009} for details on 
integration issues of IGA.} as in the traditional FEM. In words, for a macro element of order $p$,
we employ $(p+1)\times(p+1)$ Gauss points for numerical integration purposes.
Optimal quadrature rules reported in \cite{hughes_efficient_2010} use less integration points
than the standard quadrature rule (and hence results in less micro problems) are not used in this work.

We set ${{\mathbf{x}}_{{{K}_{l}}}}:=\textbf{F}\left(
\widehat{{{\mathbf{x}}_{{{K}_{l}}}}} \right)$ and consider in $\Omega $
the following sampling domains 
\begin{equation}
{{K}_{{{\delta
}_{l}}}}={{\mathbf{x}}_{{{K}_{l}}}}+\delta I
\end{equation}
where $I={{\left(
-\frac{1}{2},\frac{1}{2} \right)}^{d}}$, with $\delta \ge
\varepsilon $.  In words, the micro domain is a square centered on $\vm{x}_{K_l}$.

\subsubsection{Micro finite patch space}

On each sampling domain ${{K}_{{{\delta }_{l}}}}$, we consider a
micro triangulation ${{\mathcal{T}}_{h}}$ induced by a
geometry map which define ${{K}_{{{\delta }_{l}}}}$

\begin{equation}
{\textbf{F}_{h,{\delta _l}}}:\widehat{{{K}_{{{\delta
}_{l}}}}}\to {{K}_{{{\delta }_{l}}}},
\end{equation}

\noindent where $\widehat{{{K}_{{{\delta}_{l}}}}}$ denotes the micro parameter domain which is
$(0,1)^2$.
The partition ${{\mathcal{T}}_{h}}$ consist of elements $T$ of
diameter ${{h}_{T}}$. We denote $h:={{\max }_{K\in
{{\mathcal{T}}_{h}}}}{{h}_{T}}$.

\noindent Now, the micro finite patch space is defined as

\begin{equation}\label{microspace}
S_{{}}^{q}\left( {{K}_{{{\delta }_{l}}}},{{\mathcal{T}}_{h}}
\right)=W\left( {{K}_{{{\delta }_{l}}}} \right)\cap
\text{span}{{\left\{ {{R}_{h,ij}}\circ {{\textbf{F}}_{h,{\delta
_l}}}^{-1} \right\}}_{i,j}},
\end{equation}
where
\begin{itemize}
\item the Sobolev space $W({{K}_{{{\delta }_{l}}}})$ depends on
which kind of coupling (periodic or Dirichlet) between the macro and micro space functions
is chosen \citep{AbdulleGakuto2009}: 

\begin{equation}\label{eq:periodic_coupling}
W({{K}_{{{\delta }_{l}}}})=W_{per}^{1}({{K}_{{{\delta
}_{l}}}})=\left\{ v\in H_{per}^{1}({{K}_{{{\delta
}_{l}}}}):\int_{{{K}_{{{\delta }_{l}}}}}{vd\mathbf{x}}=0 \right\},
\end{equation}
  for periodic coupling, or
  \begin{equation}\label{eq:dirich_coupling}
W({{K}_{{{\delta }_{l}}}})=H_{0}^{1}({{K}_{{{\delta }_{l}}}})
  \end{equation}
 for Dirichlet coupling;
\item ${{\left\{ {{R}_{h,ij}} \right\}}_{i,j}}$ are the two
dimensional NURBS basis functions of order $q$ defined on the sampling domain ${{K}_{{{\delta }_{l}}}}$.
\end{itemize}
	
\noindent Next, we describe the Isogeometric Analysis Heterogeneous Multiscale Method.

\subsubsection{The IGA-HMM method}

The goal of the IGA-HMM method is to give a numerical solution
${{u}^{H}}$ of the homogenization problem \eqref{eq2}. The
variational problem is presented after introducing the macro
bilinear form and the micro problems. To facilitate the readers,
the notations in this section follow reference \citep{AbdulleGakuto2009}.

\subparagraph*{Macro-bilinear form}

We define the macro bilinear form ${{B}_{H}}$ as

\begin{equation}
{{B}_{H}}\left( {{v}^{H}},{{w}^{H}} \right):=\sum\limits_{K\in
{{\mathcal{T}}_{H}}}{\sum\limits_{l=1}^{\mathcal{L}}{\frac{{{\omega
}_{{{K}_{l}}}}}{\left| {{K}_{{{\delta }_{l}}}}
\right|}}}\int_{{{K}_{{{\delta }_{l}}}}}{{{a}^{\varepsilon
}}(\mathbf{x})\nabla v_{{{K}_{{{\delta }_{l}}}}}^{h}\cdot \nabla
w_{{{K}_{{{\delta }_{l}}}}}^{h}d\mathbf{x}},
\end{equation}

\noindent where $\mathcal{L}$ is the number of quadrature points; $v_{{{K}_{{{\delta }_{l}}}}}^{h}$ and
$w_{{{K}_{{{\delta }_{l}}}}}^{h}$ are the solutions of the
following so-called micro-problems:

\subparagraph*{Micro-problems}
Find $v_{{{K}_{{{\delta }_{l}}}}}^{h}$ (resp. $w_{{{K}_{{{\delta
}_{l}}}}}^{h}$) such that
\begin{equation}
\left( v_{{{K}_{{{\delta }_{l}}}}}^{h}-v_{lin,{{K}_{l}}}^{H}
\right)\in S_{{}}^{q}\left( {{K}_{{{\delta
}_{l}}}},{{\mathcal{T}}_{h}} \right),
\end{equation}
 and
\begin{equation}
\,\int_{{{K}_{{{\delta }_{l}}}}}{{{a}^{\varepsilon }}(\mathbf{x})\nabla
v_{{{K}_{{{\delta }_{l}}}}}^{h}\cdot {{z}^{h}}} d\vm{x}  =0,\text{
}\forall {{z}^{h}}\in {{S}^{q}}\left( {{K}_{{{\delta
}_{l}}}},{{\mathcal{T}}_{h}} \right),
\end{equation}
where $S_{{}}^{q}\left( {{K}_{{{\delta
}_{l}}}},{{\mathcal{T}}_{h}} \right)$ is the micro finite
element space defined in \eqref{microspace} and $v_{lin,{{K}_{l}}}^{H}$ is the linearization of
$v_{{{K}_{l}}}^{H}$ at the Gauss quadrature point ${{\mathbf{x}}_{{{K}_{l}}}}$ which reads

\begin{equation}
v_{lin,{{K}_{l}}}^{H}=v_{{{K}_{l}}}^{H}\left( {{\mathbf{x}}_{{{K}_{l}}}}
\right)+\nabla {{v}^{H}}\left( {{\mathbf{x}}_{{{K}_{l}}}} \right)\cdot
\left( \mathbf{x}-{{\mathbf{x}}_{{{K}_{l}}}} \right).
\end{equation}

\subparagraph*{Variational problem}
The solution ${{u}^{H}}$ of the IGA-HMM is defined by the following
variational problem: find ${{u}^{H}}\in S_{0}^{p}\left( \Omega ,{{\mathcal{T}}_{H}}
\right)$ such that

\begin{equation}\label{eq3}
{{B}_{H}}\left( {{u}^{H}},{{v}^{H}} \right)=\left\langle
f,{{v}^{H}}  \right\rangle  \equiv \int_\Omega f v^H d\vm{x} ,\text{ }\forall {{v}^{H}}\in
S_{0}^{p}\left( \Omega ,{{\mathcal{T}}_{H}} \right).
\end{equation}
More details on the implementation can be found in \citep{AbdulleShort2009}.

\subsection{A priori error estimates}

\textit{A priori} error estimates for the NURBS-based isogeometric approach have already derived in \cite{BazilevsIGAEstimates2006}, in which the authors showed that using NURBS basis functions of degree $p$ in IGA, the error convergence rates are the same as those in the standard FEM with a polynomial of the same order $p$. This fascinating result allows us to predict that the \textit{a priori} errors in the IGA-HMM follows the results derived from the FE-HMM \citep{AbdulleGakuto2009, WeinanMingZhangAnalysis2005} as shown in what follows.

Let ${{u}^{0}},{{u}^{H}}$ be the solution of problems
\eqref{eq2} and \eqref{eq3}, respectively. Under some
assumptions of regularity and periodicity on ${{a}^{\varepsilon
}}$ (\cite{AbdulleGakuto2009}, Proposition 14), using periodic
coupling \eqref{eq:periodic_coupling} and choose the size $\delta$ of sampling domain as an integer multiple of $\varepsilon$, we have the following
estimates for the solution obtained by the IGA-HMM method:
\begin{equation}\label{estimates1}
{\left\| {{u^0} - {u^H}} \right\|_{{H^1}\left( \Omega \right)}}
\leqslant C\left( {{H^p} + {{\left( {\frac{h}{\varepsilon }}
\right)}^{2q}} + \varepsilon } \right),
\end{equation}
\begin{equation}\label{estimates2}
{\left\| {{u^0} - {u^H}} \right\|_{{L^2}\left( \Omega \right)}}
\leqslant C\left( {{H^{p + 1}} + {{\left( {\frac{h}{\varepsilon
}} \right)}^{2q}} + \varepsilon } \right),
\end{equation}
{where the constant $C$ depends on the shape of the domain
$\Omega$ and the shape regularity of the mesh, but does not depend
on $\varepsilon$, the analytical solution $u^0$, and the macro/micro mesh size.}

When  the size $\varepsilon$ of the period is unknown, there is a so called cell
resonance error due to a mismatch of the sampling domain size (e.g., $\delta/\varepsilon \notin \mathbb{N}$) and to artificial boundary conditions arise, see \cite{WeinanMingZhangAnalysis2005,cellsize_2007}.

\subsection{Optimal refinement strategy for macro and micro meshes}
We denote by $N_{mac}, N_{mic}$ respectively the macro and micro
element per dimension of $\Omega $ and the sampling domains. The
corresponding macro diameter and scaled micro diameter are given
by $H=1/N_{mac}$, $h=1/N_{mic}$, respectively.

From estimates (\ref{estimates1}) and (\ref{estimates2}) we have the following optimal micro
refinement strategies \footnote{$\varepsilon$ is very small,
e.g. $\varepsilon \approx 10^{-6}$, which makes the tensor oscillate very fast, but the $u^H$ given by the HMM method can capture the homogenized solution and ``average out" these fluctuations. Also, for tensor $a^\varepsilon(x)=a(x,x/\varepsilon) $ with explicit scale separation, by collocating  slow variable at the Gauss points and choosing $\delta$ as an integer multiple of $\varepsilon$, then the $\varepsilon$ part in the right hand side of  Eqs. (\ref{estimates1}) and (\ref{estimates2}) will be ``canceled out" \citep{AbdulleShort2009}.} 

\begin{equation}\label{l2strategy}
{N_{mic}} = {N_{mac}}^{\frac{{p + 1}}{2q}}{\text{
(optimal }}{L^2}{\text{ strategy)}},
\end{equation}
\vspace{-.5cm}
\begin{equation}\label{h1strategy}
{N_{mic}} = {N_{mac}}^{\frac{p}{2q}}{\text{ (optimal }}{H^1}{\text{
strategy)}}.
\end{equation}

\subsection{Total computational cost of the IGA-HMM and the choice of macro/micro degree}

Let $N_{mac}$ denote the number of macro elements per dimension, $p,q$ the degree of macro and micro basis functions, respectively. Using the $L^2$ micro refinement strategy, we have: {Total computational cost} = {number of macro elements} $\times$ {number of micro problem per macro element} $\times$ {number of dofs per micro problem}.
Thus, the total computational cost is

\begin{equation}\label{eq:totalcost}
 N_{mac}^2 \times {(p + 1)^2} \times {(qN_{mac}^{\frac{{p + 1}}{{2q}}} + 1)^2} \sim \mathcal{O} \left( {N_{mac}^{2 + \frac{{p + 1}}{q}}} \right).
\end{equation}
\noindent From  \Eref{eq:totalcost} we see that if the degree of the micro space is fixed at $q=1$, the total cost of the IGA-HMM will increase with the degree elevation of $p$:

\begin{equation}
\text{The total cost} \sim \mathcal{O}(N_{mac}^{p+3}) ~~~ (q=1).
\end{equation}
For $p=1$ (this means linear basis functions are used in both macro and micro space), the complexity is of $\mathcal{O}(N_{mac}^{4})$ and then raise for each higher $p$, which is very inefficient. The total time to compute the solution will be very large and can be prohibited.
Therefore, instead of fixing $q=1$ at micro level, we choose $q \ge p+1 $. Then,

\begin{equation}\label{optimal_q}
\text{The total cost} \sim \mathcal{O}(N_{mac}^3).
\end{equation}
In this case, the complexity of the IGA-HMM will be hold at $\mathcal{O}(N_{mac}^3) $, for every $p$, which is acceptable  and IGA-HMM computations 
can thus be run on normal personal computers or laptops.

\begin{remark}
The formula (\ref{eq:totalcost}) also holds in FE-HMM, but in this standard approach, it is very hard to increase the basis function degree when doing implementation due to the different connection between various kinds of nodes such as center, internal, edges,...
\end{remark}


\section{Numerical examples}

In this section, we test the performance of the proposed IGA-HMM
through three thermal problems. First, a problem on a square domain with
analytical solution is analyzed with a linear basis for the micro element space. Convergence
studies are performed for various macro NURBS basis functions orders (up to order five).
Secondly, we solve a problem on a square domain that does not have an analytical solution.
As the third problem, we consider a problem defined on a curved boundary domain geometry of which can only
be exactly represented by high order NURBS. For the last two examples, high order (up to order five) NURBS elements 
are employed at both macro and micro scale which has not been done before in the framework of FE-HMM.
To the best of our knowledge, for standard FEM basis functions, in the literature, the highest degree ever used in macro space is cubic, 
and in micro space is quadratic \cite{li_efficient_2012}.

The computations are performed on a desktop computer with CPU
Intel\textsuperscript{\textregistered} core i5 2.8GHz processor.
The present method has been coded in
Matlab\textsuperscript{\textregistered} language and for
comparison purpose, the linear FE-HMM approach has also been implemented in our package.
Due to its superior performance than the Dirichlet coupling, 
periodic coupling is chosen for all examples and $\delta=\varepsilon$ is used.
A value of $\epsilon=10^{-6}$ is used for all examples. Note that the IGA-FEM or FE-HMM captures the effective solution and is thus 
independent of $\epsilon$ \citep{AbdulleShort2009}. $C^0$ NURBS are used in the micro spaces unless otherwise stated.
We recall  that $C^0$ NURBS are similar to high order finite elements ($C^0$ continuity across element boundaries), refer to Fig. \ref{fig:CkC0}.

\subsection{\textbf{IGA-HMM on domain with straight boundary}}
\label{sec:problem1}

We consider the following two-scale problem \cite{AbdulleGakuto2009}

\begin{equation}
\begin{gathered}
  \left\{ \begin{gathered}
- \nabla \cdot\left( {a\left( {\frac{\mathbf{x}}{\varepsilon}}
\right)\nabla {u^\varepsilon}(\mathbf{x})} \right) = 1{\text{ in }}\Omega
= {\left( {0,1} \right)^2}, \\
{u^\varepsilon}(\mathbf{x}) = 0{\text{ on }}\partial {\Omega _D}: = \{
{x_1} = 0\} \cup \left\{ {{x_2} = 1} \right\} , \\
\mathbf{n}\cdot\left( {a\left( {\frac{\mathbf{x}}{\varepsilon}} \right)\nabla
{u^\varepsilon}(\mathbf{x})} \right) = 0{\text{ on }}\partial {\Omega
_N}: = \partial \Omega \backslash \partial {\Omega _D} ,
\end{gathered}  \right.
\end{gathered}
\end{equation}
where $a\left( {\frac{\mathbf{x}}{\varepsilon}} \right) = (\cos \left( {2\pi \frac{{{x_1}}}{\varepsilon}} \right) + 2)\mathbf{I_2},  \mathbf{x}=(x_1,x_2), \text{and }\mathbf{I_2} \text{ is the identity matrix}. $  

The corresponding homogenized tensor and
solution can be computed analytically \\
\[u^0(\vm{x})=\frac{-x_1^2+x_1}{2\sqrt{3}},\;\;\; {a^0} = \left[ {\begin{array}{*{20}{c}}
  {{{\left( {\int_0^1 {\frac{{dt}}{{k(t)}}} } \right)}^{ - 1}}}&0 \\
  0&2
\end{array}} \right]=\left[ {\begin{array}{*{20}{c}}
  \sqrt{3}&0 \\
  0&2
\end{array}} \right],\]
where $k(t) = \cos 2\pi t + 2$. The analytical solution $u_0$ is plotted in Fig. \ref{fig:square_solution}.

We test the performance of IGA-HMM by considering the
${{L}^{2}}$ error and ${{H}^{1}}$ (energy) errors. 

In this example, we choose the degree $q$ of the micro basis functions equal to 1. Hence, the micro element becomes $Q4$ as in FE-HMM. 
We study the convergence of the error measured in $H^1$ and $L^2$ norm with different NURBS basis order $p$ in the macroscopic solver.
Set $q=1$ from (\ref{l2strategy}), (\ref{h1strategy}) we have

\begin{center}
$\begin{gathered}
{N_{mic}} = {N_{mac}}^{\frac{{p + 1}}{2}}{\text{
(}}{L^2}{\text{ strategy)}} ,\hfill \\
{N_{mic}} = {N_{mac}}^{\frac{p}{2}}{\text{ (}}{H^1}{\text{ strategy)}}. \hfill \\
\end{gathered} $
\end{center}

\begin{figure}[H]
\begin{center}
\subfigure{\includegraphics[trim= 2cm 1cm 2cm 1cm,
angle=0,width=0.45\columnwidth]{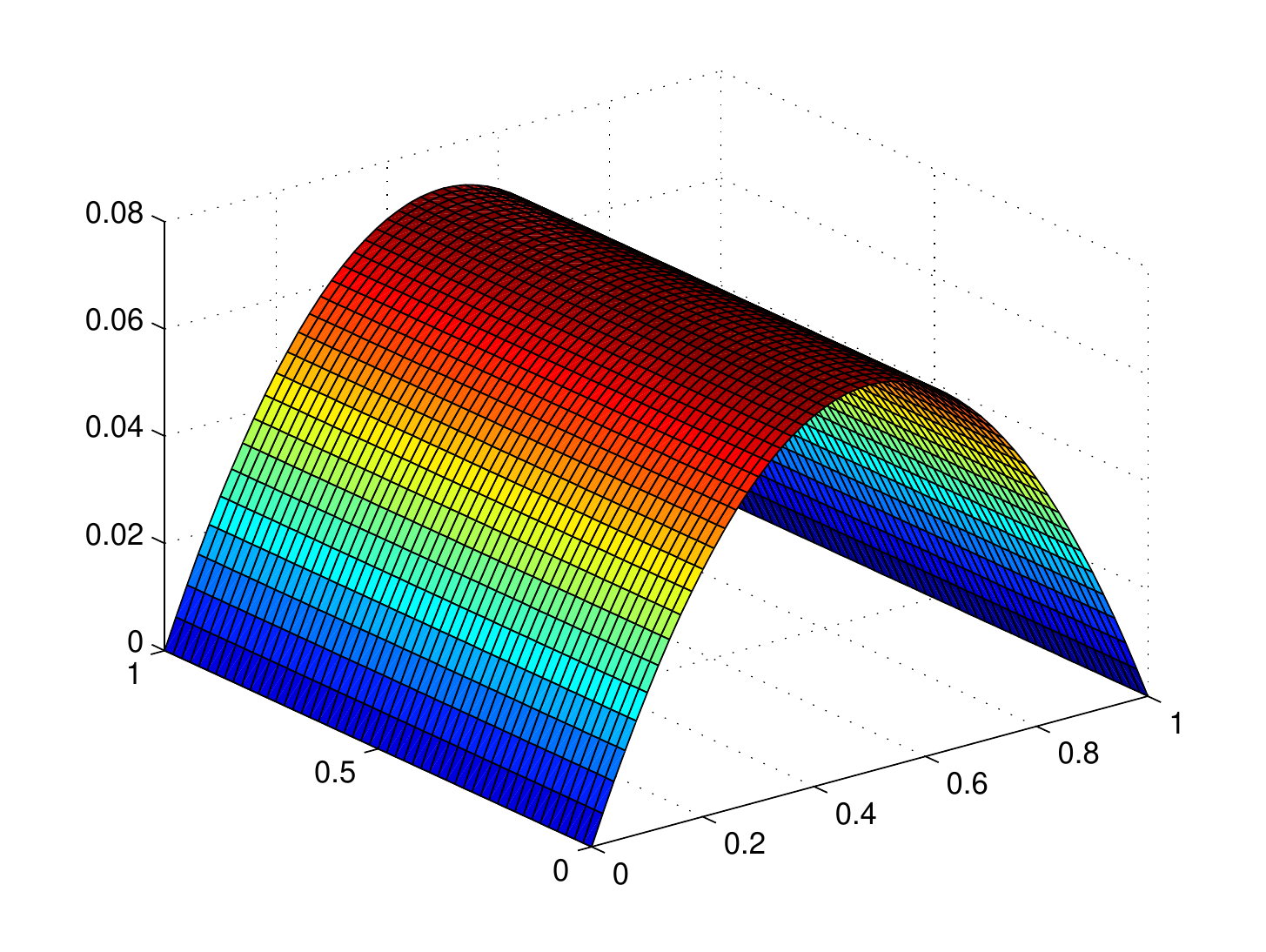}}
\end{center}
\caption{Homogenized solution of the problem ~\ref{sec:problem1}.}
\label{fig:square_solution}
\end{figure}

\subsubsection{IGA-HMM: convergence rate test}

First, we study the \textit{$L^2$ micro refinement
strategy}. From \Fref{fig:square_errL2_L2strategy} we see that
for linear, quadratic, cubic, and quartic NURBS, the
convergence rates of the errors in $L^2$ norm when using the $L^2$
micro refinement strategy have the slopes of $2,3,4,5$,
respectively. These match well the predicted theory. The details
of the errors are given in Tab. \ref{tab:square_errL2_L2strategy}.

\begin{table}[H]
\centering
\begin{threeparttable}
\begin{tabularx}{\textwidth}{XCCC}
\toprule
&\multicolumn{3}{c}{Mesh}\\ 
\cline{2-4}
Method & 2$\times$2 & 4$\times$4&8$\times$8 \\
\midrule
FE-HMM\tnote{a}			&0.28950&	0.08360	&0.02100\\
p=2 (IGA-HMM)&0.06092	&8.36e-3	&1.03e-3\\
p=3 (IGA-HMM)&0.03452	&2.13e-3	&1.34e-4\\
p=4 (IGA-HMM)&0.01443	&5.34e-4	&1.66E-05\\
\bottomrule
\end{tabularx}
\begin{tablenotes}
       \item[a] In this work, the standard $Q4$ elements are used
for FE-HMM in both macro and micro finite element space.
     \end{tablenotes}
  \end{threeparttable}
\caption{$L^2$ error of problem \ref{sec:problem1}, using the $L^2$ micro refinement strategy.}
\label{tab:square_errL2_L2strategy}
\end{table}

\begin{figure}[H]
\begin{center}
\subfigure{\includegraphics[trim= 2cm 1cm 2cm 1cm,
angle=0,width=.45\columnwidth]{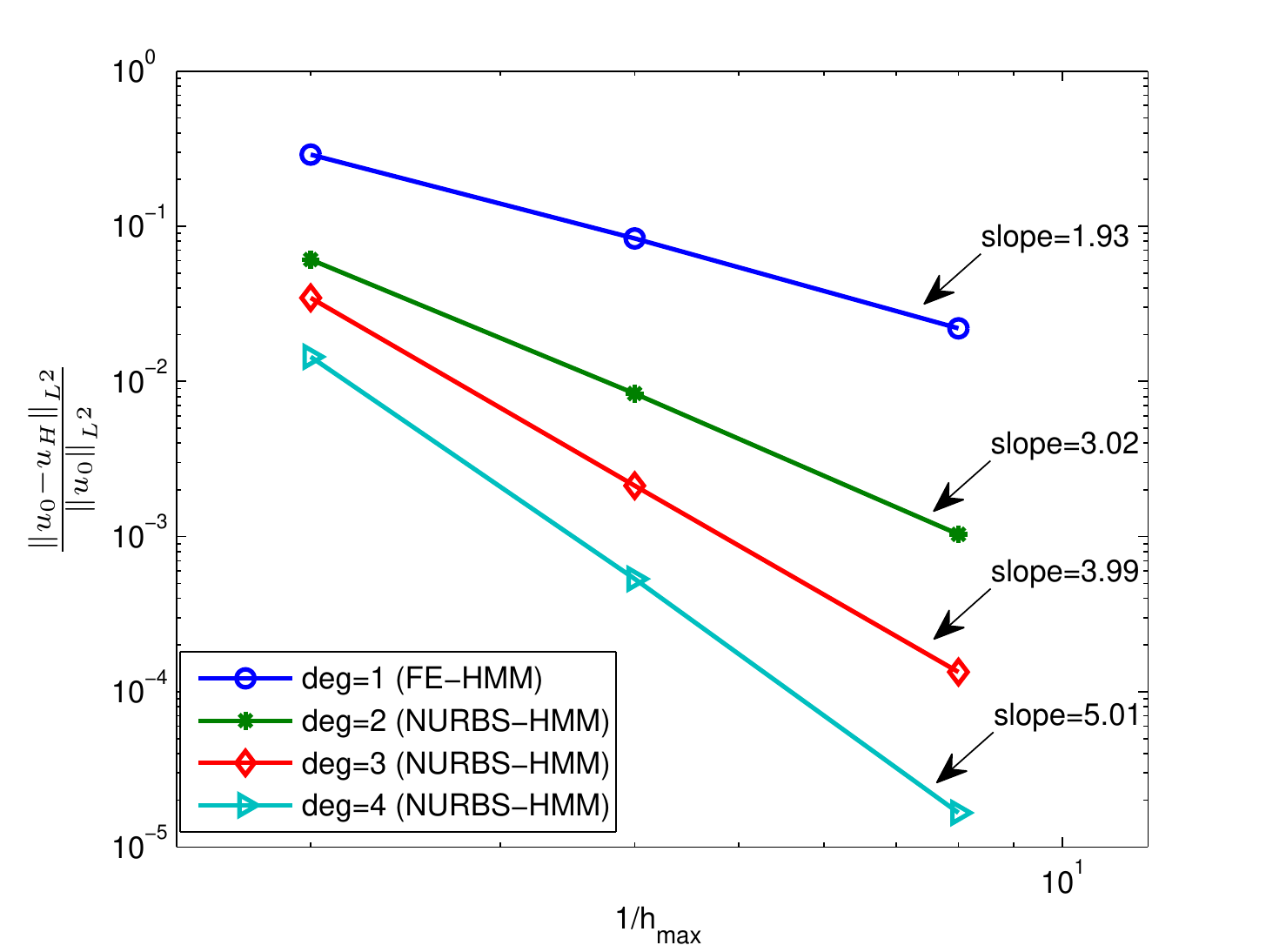}}
\end{center}
\caption{$L^2$ error of problem \ref{sec:problem1}, using the $L^2$
micro refinement strategy}
\label{fig:square_errL2_L2strategy}
\end{figure}

For the $H^1$ error when using the $L^2$ micro refinement strategy,
from \Fref{fig:square_errH1_L2strategy}, we see that the
convergence rate is slightly higher than the theoretical value. This is due
{\color{black}{to the problem dependence. In the later examples, this phenomenon does not occur.}} We list in Tab. \ref{tab:square_errH1_L2strategy} the detailed errors.

\begin{table}[h]
\centering
\begin{tabularx}{\textwidth}{XCCC}
\toprule
&\multicolumn{3}{c}{Mesh}\\ 
\cline{2-4}
{Method}& 2$\times$2 & 4$\times$4&8$\times$8\\
\hlinewd{0.8pt}
FE-HMM   & 0.50559	 &  0.25168	 & 0.12527\\
p=2 (IGA-HMM) & 0.06092	& 8.36e-3	& 1.03e-3\\
p=3 (IGA-HMM) & 0.03435	& 2.13e-3	& 1.34e-4\\
p=4 (IGA-HMM) & 0.01443	& 5.34e-4	& 1.66E-05\\
\bottomrule
\end{tabularx}
\caption{$H^1$ error of problem \ref{sec:problem1} using the $L^2$ micro refinement strategy.}
\label{tab:square_errH1_L2strategy}
\end{table}

\begin{figure}[H]
\begin{center}
\subfigure{\includegraphics[trim= 2cm 1cm 2cm 1cm,
angle=0,width=0.45\columnwidth]{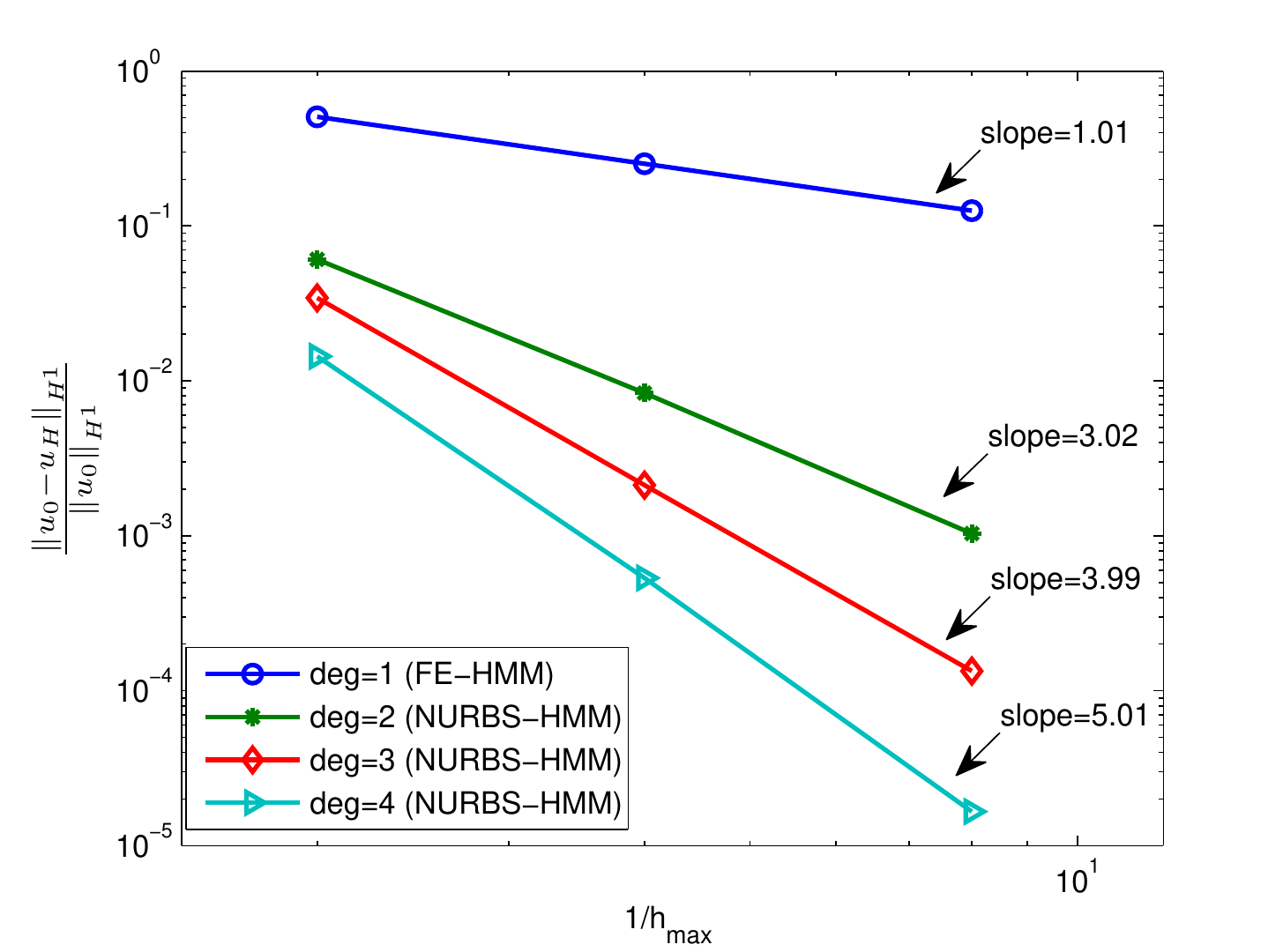}}
\end{center}
\caption{$H^1$ error of the problem \ref{sec:problem1}, using the $L^2$
micro refinement strategy}
\label{fig:square_errH1_L2strategy}
\end{figure}

We next continue to test the problem for the \textit{$H^1$ micro
refinement strategy}. \Fref{fig:square_errH1_H1strategy}  and Tab. \ref{tab:square_errH1_H1strategy} 
demonstrate that for linear,
quadratic, cubic, quartic, and quintic NURBS, the convergence rates
of errors in $H^1$ norm are $1,2,3,4,5$ and in $L^2$ norm (see Fig. \ref{fig:square_errL2_H1strategy} and Tab. \ref{tab:square_errL2_H1strategy}) are
$1,2,3,4,5$, respectively. It is clear that the results 
match the theoretical ones.
\begin{table}[!ht]
\centering
\begin{tabularx}{\textwidth}{lCCCCCCC}
\toprule
&\multicolumn{6}{c}{Mesh}\\ \cline{2-7}
Method & 2$\times$2 & 4$\times$4&
8$\times$8 & 16$\times$16 & 32$\times$32 & 64$\times$64 \\
\midrule
FE-HMM & 0.50559 & 0.27581 & 0.13888 &0.06725 & 0.03443 &
0.01770\\
p=2 (IGA-HMM) & 0.07647 & 0.02964 & 8.31e-3 & 2.13e-3 &
5.34e-4&-- \\
p=3 (IGA-HMM)& 0.06439	& 8.34e-3	& 1.03e-3 & 1.34e-4 &-- &-- \\
p=4 (IGA-HMM)& 0.03108	& 2.13e-3	&  1.34e-4 &-- &-- &-- \\
p=5 (IGA-HMM)& 0.01414	& 5.34e-4	& 1.66e-05 &-- &-- & --\\
\bottomrule
\end{tabularx}
\caption{$H^1$ error of the problem \ref{sec:problem1}, using the $H^1$
micro refinement strategy.}
\label{tab:square_errH1_H1strategy}
\end{table}
\begin{figure}[H]
\begin{center}
\includegraphics[trim= 2cm 1cm 2cm 1cm,
angle=0,width=0.45\columnwidth]{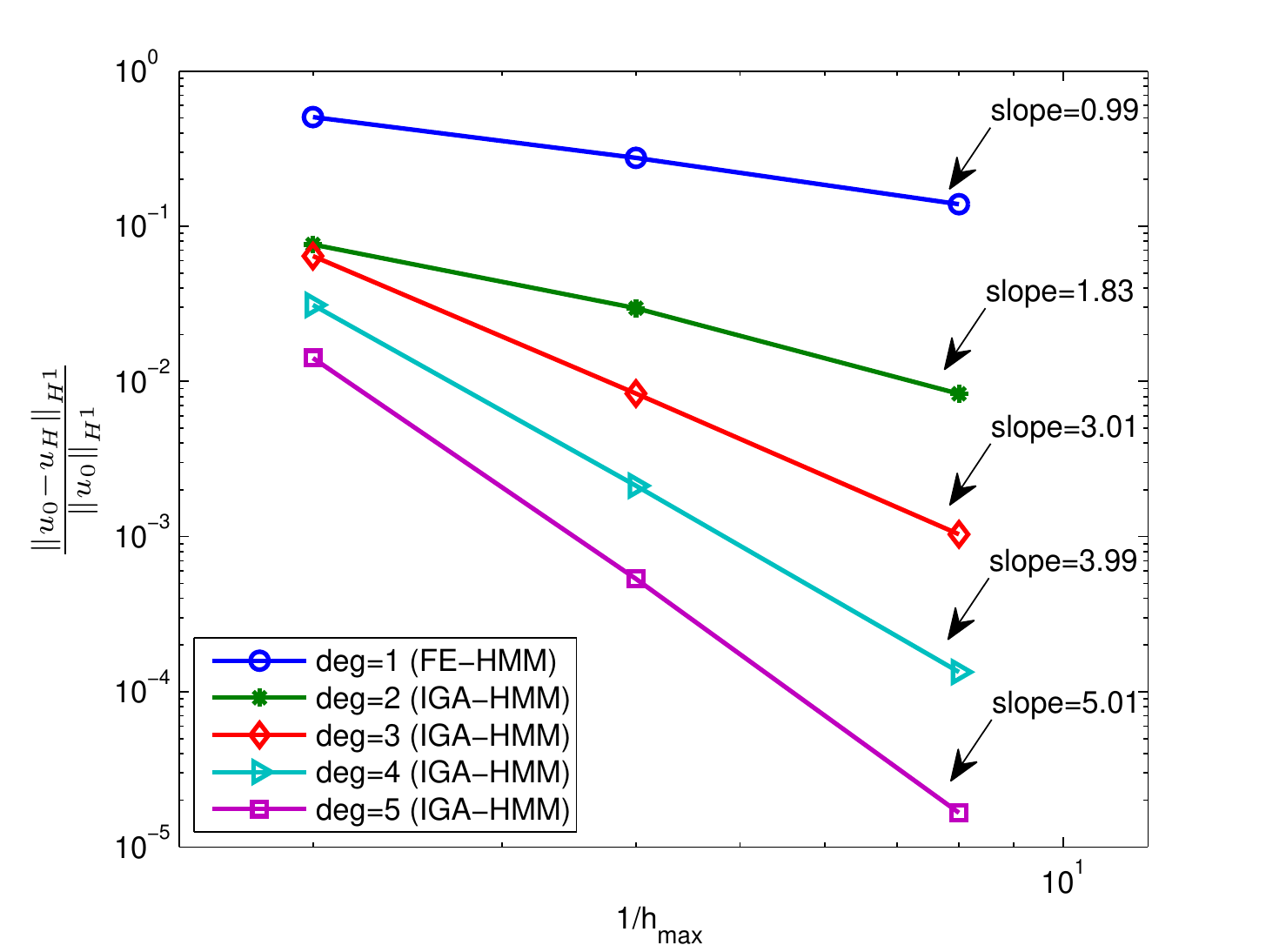}
\end{center}
\caption{$H^1$ error of problem \ref{sec:problem1}, using the $H^1$
micro refinement strategy.}
\label{fig:square_errH1_H1strategy}
\end{figure}
%
\begin{table}[H]
\centering
\begin{tabularx}{\textwidth}{lCCCCCCC}
\toprule
&\multicolumn{6}{c}{Mesh}\\ \cline{2-7}
Method & 2$\times$ 2 & 4$\times$4& 8$\times$8 & 16$\times$16 & 32$\times$32 & 64$\times$64\\
\midrule
FE-HMM & 0.28950 & 0.16737 & 0.07342&0.02809 &0.015260
&8.51e-3\\
p=2 (IGA-HMM)& 0.07646	& 0.02968	& 8.31e-3&2.13e-3	&5.34e-4&--\\p=3 (IGA-HMM)& 0.06471	& 8.35e-3	& 1.03e-3&0.00013&--&--	\\
p=4 (IGA-HMM)& 0.03120	& 2.13e-3	& 1.34e-4&--&--&--\\
p=5 (IGA-HMM)& 0.01413	& 5.34e-4	& 1.66e-05&--&--&--\\
\bottomrule
\end{tabularx}
\caption{$L^2$ error of the thermal square problem, using the $H^1$
micro refinement strategy.}
\label{tab:square_errL2_H1strategy}
\end{table}

\begin{figure}[H]
\begin{center}
\subfigure{\includegraphics[trim= 2cm 1cm 2cm 1cm,
angle=0,width=.45\columnwidth]{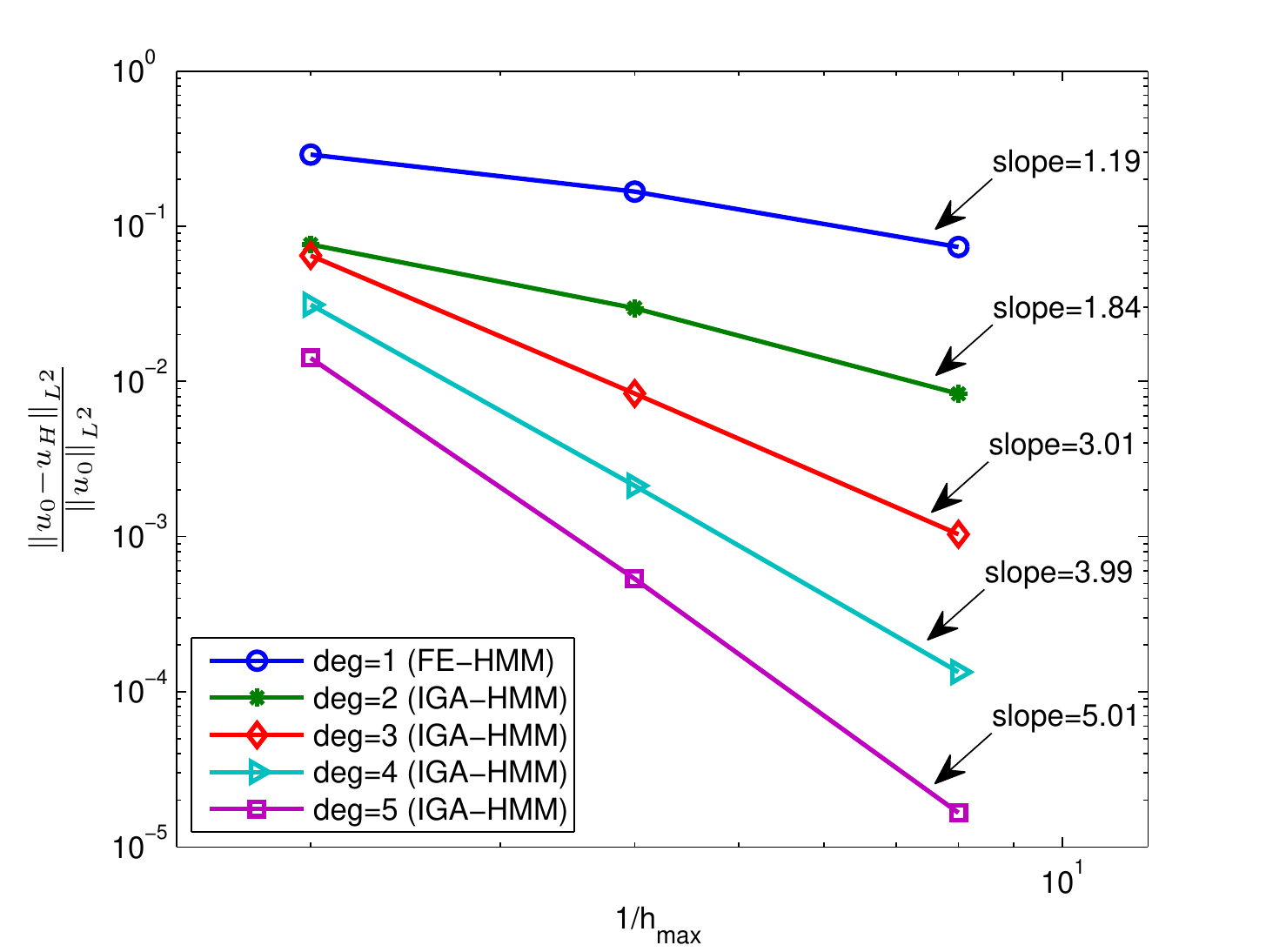}}
\end{center}
\caption{$L^2$ error of problem \ref{sec:problem1}, using the $H^1$ micro refinement strategy.}
\label{fig:square_errL2_H1strategy}
\end{figure}

\subsubsection{IGA-HMM vs FE-HMM: performance comparison}

In this section we compare the CPU-time of the IGA-HMM and that of
the FE-HMM needed to obtain solutions with various 
accuracies: $10^{-2}, 2 \times 10^{-3}, 5\times 10^{-4},10^{-4}, \text{and } 10^{-5}$ in $H^1 $ norm of the error, using the $H^1$ micro
refinement strategy. The results are presented in
Tab. \ref{tab:compare_target}. 
We can see that for the same accuracy, the IGA-HMM
outperforms the (linear) FE-HMM in terms of CPU-time needed to compute
the solution. For example, to obtain the $H^1$ error less than $10^{-2}$,
the FE-HMM needs about 1300 second, while the IGA-HMM needs only
approximately 5 seconds and 2.5 seconds for quadratic and cubic
NURBS, respectively. 
Similarly, to obtain the error in $H^1$ norm approximately
$2\times 10^{-3}$, the FE-HMM needs more than 80000 seconds, while
the IGA-HMM needs only approximately 83 seconds and 16 seconds
for quadratic and quartic NURBS, respectively. 
We can see that the CPU time of the IGA-HMM is significantly reduced
in comparison with the FE-HMM thanks to the high order NURBS basis functions
employed at the macroscale which allows very coarse macro meshes to be used.
It should be noted that our Matlab implementation was not optimized and hence performance can be further
improved.

\begin{table}[!ht]
\centering
\begin{tabularx}{\textwidth}{lXXX}
\toprule
 method &mesh&$H^1$ err & CPU(s)\\
\midrule
FE-HMM    &$128\times128$&8.67e-3	    &1359\\
p=2 (IGA-HMM)&$8\times8$& 8.31e-3	& 4.9\\
p=3 (IGA-HMM)&$4\times4$&  8.35e-3  & 2.5\\
\hline
FE-HMM    &$512\times512$&2.17e-3	    &$>$80000 \\
p=2 (IGA-HMM)&$16\times16$& 2.13e-3	& 83\\
p=4 (IGA-HMM)&$4\times4$&  2.13e-3  & 16\\
\hline
FE-HMM    &$2048\times2048$&5.42e-3	    &$>$5e6 (est)\\
p=2 (IGA-HMM)&$32\times32$& 5.34e-4	& 1433\\
p=5 (IGA-HMM)&$4\times4$&  5.34e-4  & 102\\
\hline
FE-HMM    &$8192\times8192$&1.35e-4	    &$>$3e8 (est)\\
p=3 (IGA-HMM)&$16\times16$& 1.34e-4	& 3064\\
p=6 (IGA-HMM)&$4\times4$&  1.34e-4  & 673\\
\hline
FE-HMM    &$65536\times 65536$&1.69e-5	    &$>$2e10(est)\\
p=5 (IGA-HMM)&$8\times8$& 1.66e-5	& 6.06e4\\
\bottomrule
\end{tabularx}
\caption{{\color{black}Performance comparison between the FE-HMM and the IGA-HMM (both using linear basis functions in micro space).}}
\label{tab:compare_target}
\end{table}

\subsection{\textbf{A high-order approximation of IGA-HMM in both macro and micro levels
patch spaces}}\label{sec:problem2}

In the previous example, we have already seen that the IGA-HMM
performs much better than FE-HMM when using high order NURBS ($p\ge 2$) in macro patch space and NURBS with degree $q=1$ in micro
space. In this section, we continue to test the IGA-HMM with high order NURBS in both micro and macro space.
To this end, we consider the following two-scale problem:

\begin{equation}
\begin{gathered}
  \left\{ \begin{gathered}
- \nabla \cdot\left( {a\left(\mathbf{x}, {\frac{\mathbf{x}}{\varepsilon}}
\right)\nabla {u^\varepsilon}(\mathbf{x})} \right) = 1{\text{ in }}\Omega
= {\left( {0,1} \right)^2}, \\
{u^\varepsilon}(\mathbf{x}) = 0{\text{ on }}\partial {\Omega _D}: = \{
{x_1} = 0\} \cup \left\{ {{x_2} = 1} \right\} , \\
\mathbf{n}\cdot\left( {a\left(\mathbf{x}, {\frac{\mathbf{x}}{\varepsilon}} \right)\nabla
{u^\varepsilon}(\mathbf{x})} \right) = 0{\text{ on }}\partial {\Omega
_N}: = \partial \Omega \backslash \partial {\Omega _D} ,
\end{gathered}  \right.
\end{gathered}
\end{equation}
where the conductivity tensor is given by \cite{abdulle_reduced_2012}:

\begin{equation}
a\left(\mathbf{x}, {\frac{\mathbf{x}}{\varepsilon}}
\right)=\left[ {\begin{array}{*{20}{c}}
{x_1^2 + 0.2 + \left( {{x_2} + 1} \right)\left( {\sin \left( {2\pi \frac{{{x_1}}}{\varepsilon }} \right) + 2} \right)}&0\\
0&{x_2^2 + 0.05 + \left( {{x_1}{x_2} + 1} \right)\left( {\sin \left( {2\pi \frac{{{x_2}}}{\varepsilon }} \right) + 2} \right)}
\end{array}} \right].
\end{equation}
The corresponding homogenized tensor is given by 
\begin{equation}
{a^0}(\mathbf{x}) = \left[ {\begin{array}{*{20}{c}}
{{{\left( {\int_0^1 {\frac{{d{y_1}}}{{x_1^2 + 0.2 + \left( {{x_2} + 1} \right)\left( {\sin \left( {2\pi {y_1}} \right) + 2} \right)}}} } \right)}^{ - 1}}}&0\\
0&{{{\left( {\int_0^1 {\frac{{d{y_2}}}{{x_2^2 + 0.05 + \left( {{x_1}{x_2} + 1} \right)\left( {\sin \left( {2\pi {y_2}} \right) + 2} \right)}}} } \right)}^{ - 1}}}
\end{array}} \right].
\end{equation}
This problem does not have an explicit analytical solution, thus for a reference solution, we solve the homogenized problem with the homogenized tensor $a^0$ on a fine mesh of $500\times 500$ {\color{black}quintic elements using the standard NURBS-based finite element method, see \Fref{fig:abdullesquare_solution}.
The relative errors are computed using the following formula:}
\[er{r_{{H^1}}} = \frac{{{{\left\| {{u_H} - {u_{ref}}} \right\|}_{{H^1}(\Omega )}}}}{{{{\left\| {{u_{ref}}} \right\|}_{{H^1}(\Omega )}}}}, ~~~ er{r_{{L^2}}} = \frac{{{{\left\| {{u_H} - {u_{ref}}} \right\|}_{{L^2}(\Omega )}}}}{{{{\left\| {{u_{ref}}} \right\|}_{{H^1}(\Omega )}}}}.\]

For the macro space we use NURBS of degree
ranging from 1 to 3, and for the micro space, from the choice $q \ge p+1 $ as discussed in Eq. (\ref{optimal_q}), we use the NURBS function of the degree fixed at 5 ($q$ should not be too large, because when the micro mesh is too coarse, the results will not stable). The results are presented in Tabs. \ref{tab:abdulle_errH1_L2strategy} and \ref{tab:abdulle_errL2_L2strategy}.

\begin{figure}[!ht]
\centering
\includegraphics[angle=0,width=.6\columnwidth]{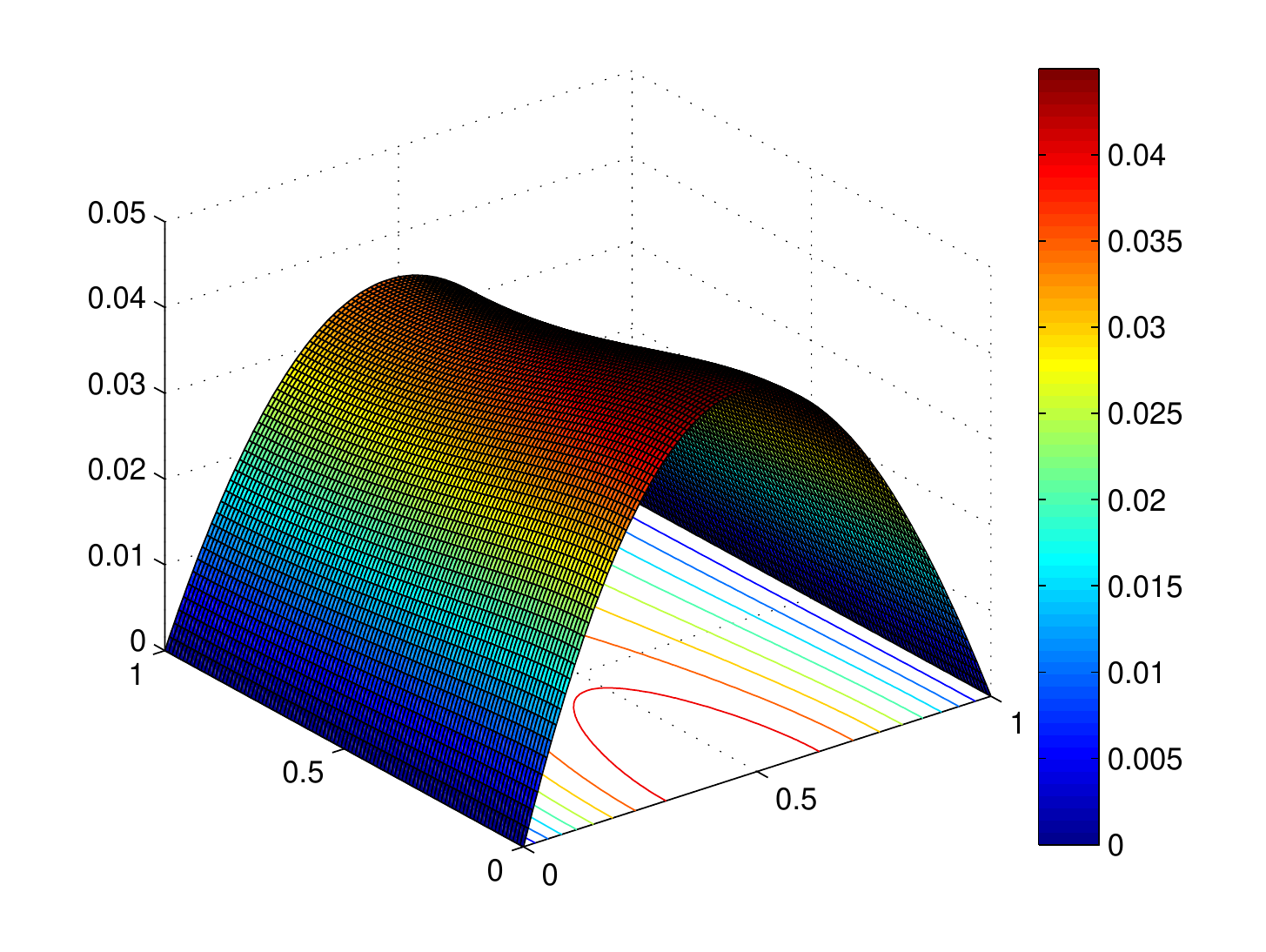}
\caption {{\color{black}Homogenized solution of the problem \ref{sec:problem2}.}
\label{fig:abdullesquare_solution}}
\end{figure}

\begin{table}[H]
\centering
\begin{tabularx}{\textwidth}{lXXXXX}
\toprule
&\multicolumn{4}{c}{Mesh}\\ \cline{2-6}
{Degree}      & $4\times4$  & $8\times8$ & $16\times16$ & $32\times32$  & $64\times64$ \\ 
\midrule
FE-HMM & 0.24884&	0.12443&	0.06223	&0.03112 & 0.01557\\
p=2 &6.97 e-3	&1.71 e-3	&4.26	e-4&1.06 e-4 & 2.66e-5\\
p=3 &6.73 e-4	&9.01 e-5	&1.20	 e-5 &1.59 e-6 & 2.07e-7\\
\bottomrule
\end{tabularx}
\caption{$H^1$ error of problem \ref{sec:problem2}, using the $L^2$ micro refinement strategy}
\label{tab:abdulle_errH1_L2strategy}
\end{table}

\begin{table}[H]
\centering
\begin{tabularx}{\textwidth}{lXXXXX}
\toprule
&\multicolumn{4}{c}{Mesh}\\ \cline{2-6}
{Degree}      & $4\times4$  & $8\times8$ & 1$6\times16$ & $32\times32$ & $64\times64$ \\ 
\midrule
FE-HMM &0.06063&	0.01514	&0.00379&	0.00095 & 0.00024\\
p=2 &7.66 e-4	&8.95 e-5	&1.09 e-5	&1.37 e-6 & 1.72e-7\\
p=3& 6.51 e-5	&4.10 e-6	&2.70e-7	&3.24 e-8 & 1.11e-9\\
\bottomrule
\end{tabularx}
\caption{$L^2$ error of problem \ref{sec:problem2}, using the $L^2$ micro refinement strategy.}
\label{tab:abdulle_errL2_L2strategy}
\end{table}

\begin{figure}[H]
\begin{center}
\subfigure{\includegraphics[trim= 2cm 1cm 2cm 1cm,
angle=0,width=.45\columnwidth]{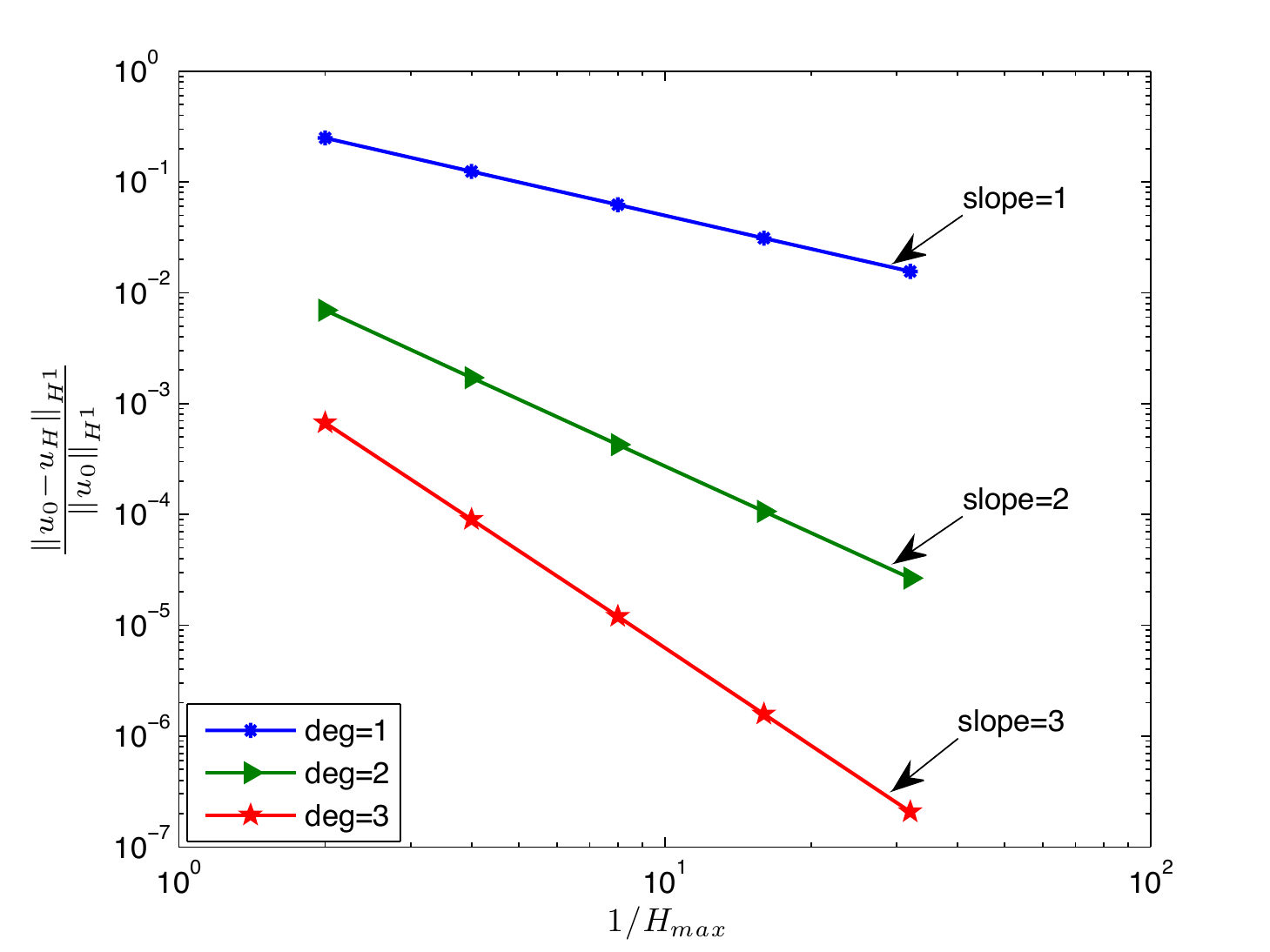}}
\end{center}
\caption{{\color{black}{$H^1$ error of problem \ref{sec:problem2}, quintic NURBS in micro space and $L^2$ refinement strategy are used.}}}
\label{fig:abdulle_errH1_L2strategy}
\end{figure}

\begin{figure}[H]
\begin{center}
\subfigure{\includegraphics[trim= 2cm 1cm 2cm 1cm,
angle=0,width=.45\columnwidth]{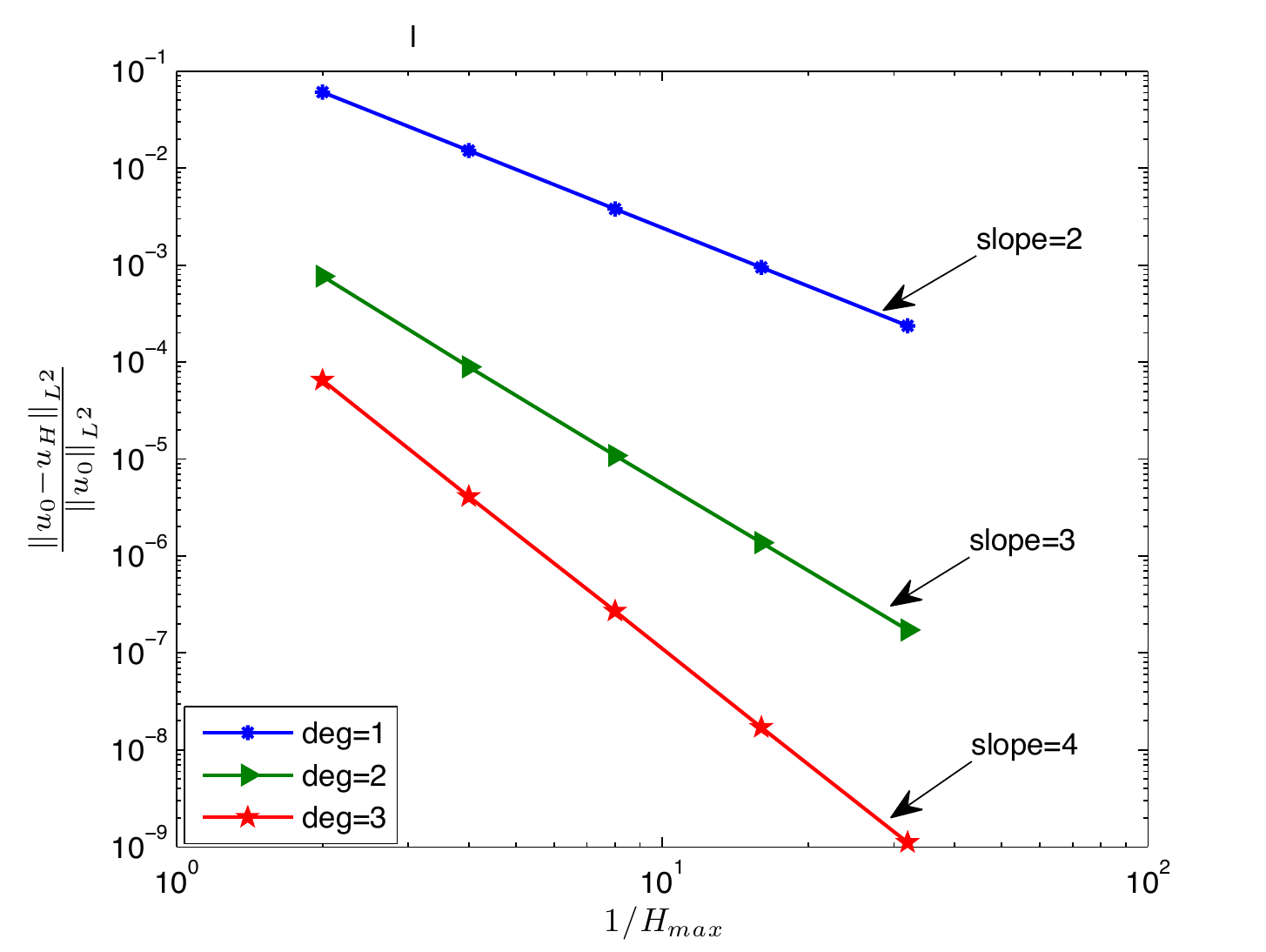}}
\end{center}
\caption{{\color{black}{$L^2$ error of problem \ref{sec:problem2}, quintic NURBS in micro space and $L^2$ refinement strategy are used.}}}
\label{fig:abdulle_errL2_L2strategy}
\end{figure}

Figures \ref{fig:abdulle_errH1_L2strategy} and \ref{fig:abdulle_errL2_L2strategy} shows that the convergence rate of the IGA-HMM is optimal and matches the theory suggested by FE-HMM: $1,2,3$ in $H^1$ error and $2,3,4$ in $L^2$ error for linear, quadratic and cubic NURBS in the macro space, respectively. From Tabs.\ref{tab:abdulle_errH1_L2strategy} and \ref{tab:abdulle_errL2_L2strategy}, we can see that the IGA-HMM is able to obtain a very high accuracy, up to $2.07e-7$ and $1.11 e-9$ in $H^1$ error and $L^2$ error, respectively. Such accuracy is almost prohibited in FE-HMM, as pointed out in the previous example. Also in comparison with the results obtained in the previous example, we can see that the high order approach in both micro and macro spaces is more effective than when we only use high order in macro and linear basis functions in micro spaces, which is only up to $10^{-6}$ in $L^2$ error.  


\subsection{\textbf{IGA-HMM applied for a curved boundary domain}}
\label{sec:problem3}

In this problem, we test the convergence of the IGA-HMM method
on a curved boundary domain, where traditional FEM methods
cannot exactly represent the boundary. Let $\Omega$ be the
quarter circular annulus of internal radius $r=1$ and external
radius $R=2$, which lies in the first quadrant of the Cartesian
plane, see \Fref{fig:quartercylinder_coarse}. Homogeneous
Dirichlet boundary condition is imposed on the whole boundary.
Let the heat source be

\begin{equation}
f=-{\frac {2\,{x_1}^{3}x_2 \left(18\,\sqrt {{x_1}^{2}+{x_2}^{2}}+4\,\sqrt {3}-
24 \right) +2\,x_1{x_2}^{3} \left( 9\,\sqrt {3}\sqrt {{x_1}^{2}+{x_2}^{2}}-12
\,\sqrt {3}+8 \right)}{ \left( {x_1}^{2}+{x_2}^{2} \right) ^{3}}}.
\end{equation}
We consider the following problem:

\begin{equation}
\begin{gathered}
  \left\{ \begin{gathered}
- \nabla \cdot\left( {a\left( {\frac{\mathbf{x}}{\varepsilon}}
\right)\nabla {u^\varepsilon}(\mathbf{x})} \right) = f{\text{ in }}\Omega
, \\
  {u^\varepsilon}(\mathbf{x}) = 0{\text{  on }}\partial {\Omega} \\
\end{gathered}  \right.
\end{gathered}
\end{equation}
where $a\left( {\frac{\mathbf{x}}{\varepsilon}} \right) = (\cos \left( {2\pi \frac{{{x_1}}}{\varepsilon}} \right) + 2)\mathbf{I_2}, \mathbf{x}=(x_1,x_2), \text{and }\mathbf{I_2} \text{ is the identity matrix}.$ 

The analytical homogenized solution $u^0$ is given by (see also \Fref{fig:quartercylinder_solution})

\begin{equation}
u^0=2x_1x_2+{\frac {x_1x_2 \left( 4-6\sqrt {{x_1}^{2}+{x_2}^{2}}
\right)}{{x_1}^{2}+{x_2}^{2}}}.
\end{equation}
First, we consider the case with linear basis function at micro space and for three macro meshes
as shown in \Fref{fig:quartercylinder_meshes}. The geometry data for the coarsest mesh (single quadratic element) is given in
Tab. \ref{tab:annulus_data}. Other meshes of different densities and basis orders are constructed seamlessly from this coarsest mesh. 
Tables~\ref{tab:annulus_errH1_H1strategy},~\ref{tab:annulus_errL2_H1strategy}
and ~\Fref{fig:annulus_errH1_H1strategy}, ~\ref{fig:annulus_errL2_H1strategy} demonstrate that the convergence rate of the error is optimal and matches the theory.

\begin{figure}[H]
\centering
\includegraphics[trim= 2cm 1cm 0cm 0cm ,
angle=0,width=.45\columnwidth]{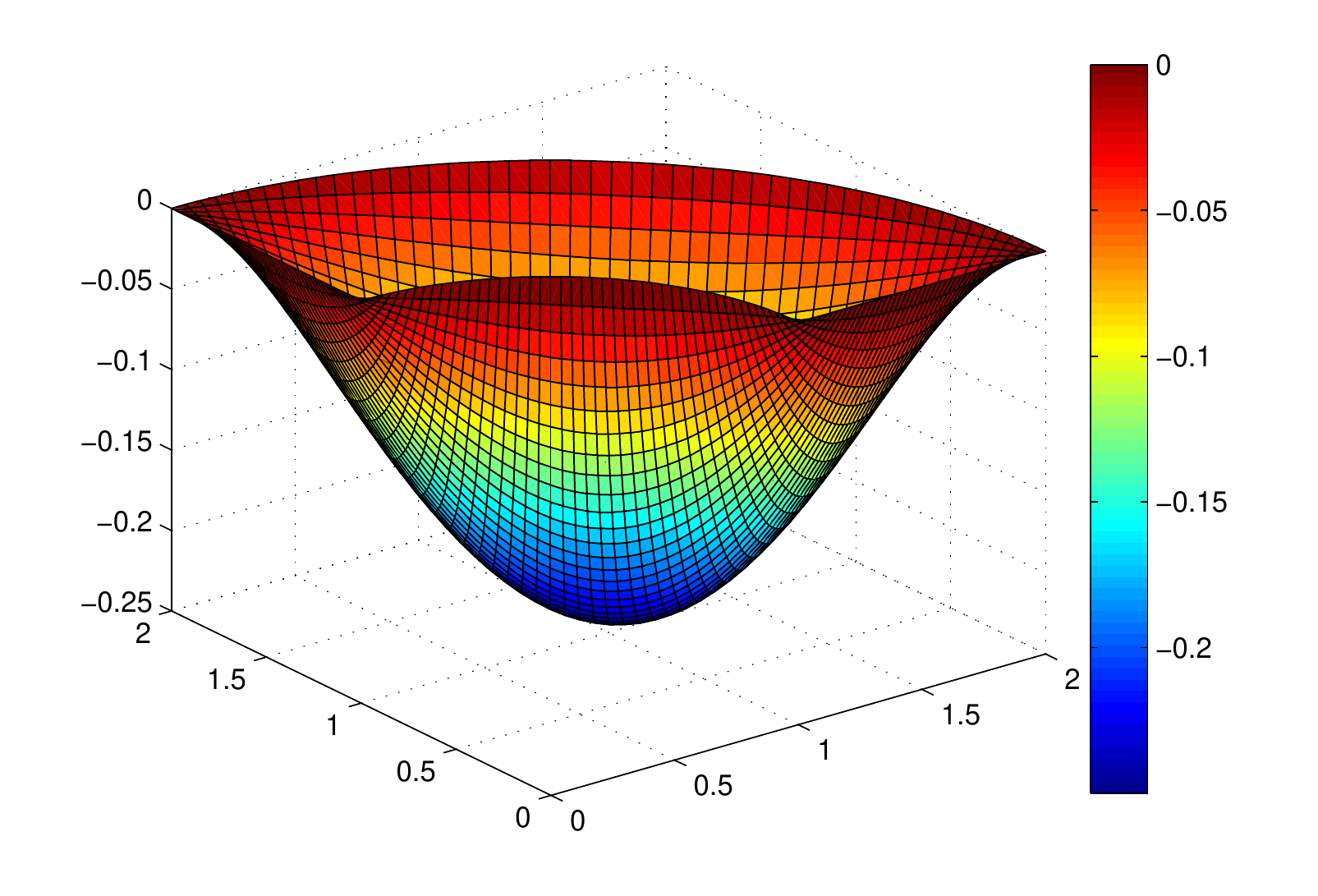}
\includegraphics[trim= 2cm 1cm 0cm 0cm,
angle=0,width=.45\columnwidth]{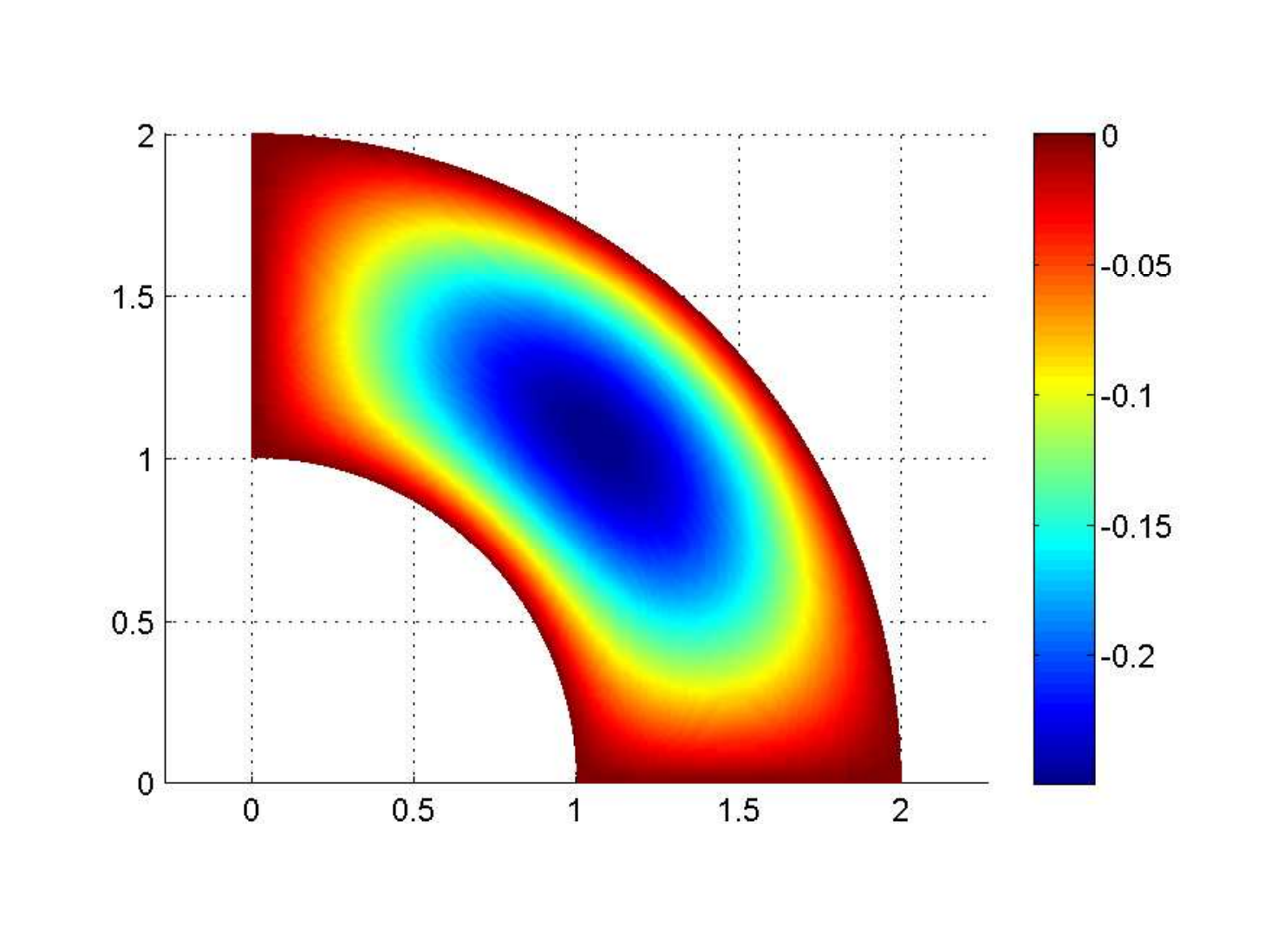}
\caption {Homogenized solution of the problem~~\ref{sec:problem3} (left: 3D and right: contour plot).}
\label{fig:quartercylinder_solution}
\end{figure}

\begin{figure}[!ht]
\centering
\includegraphics[trim=0cm 1cm 2cm 1cm,width=.43\columnwidth]{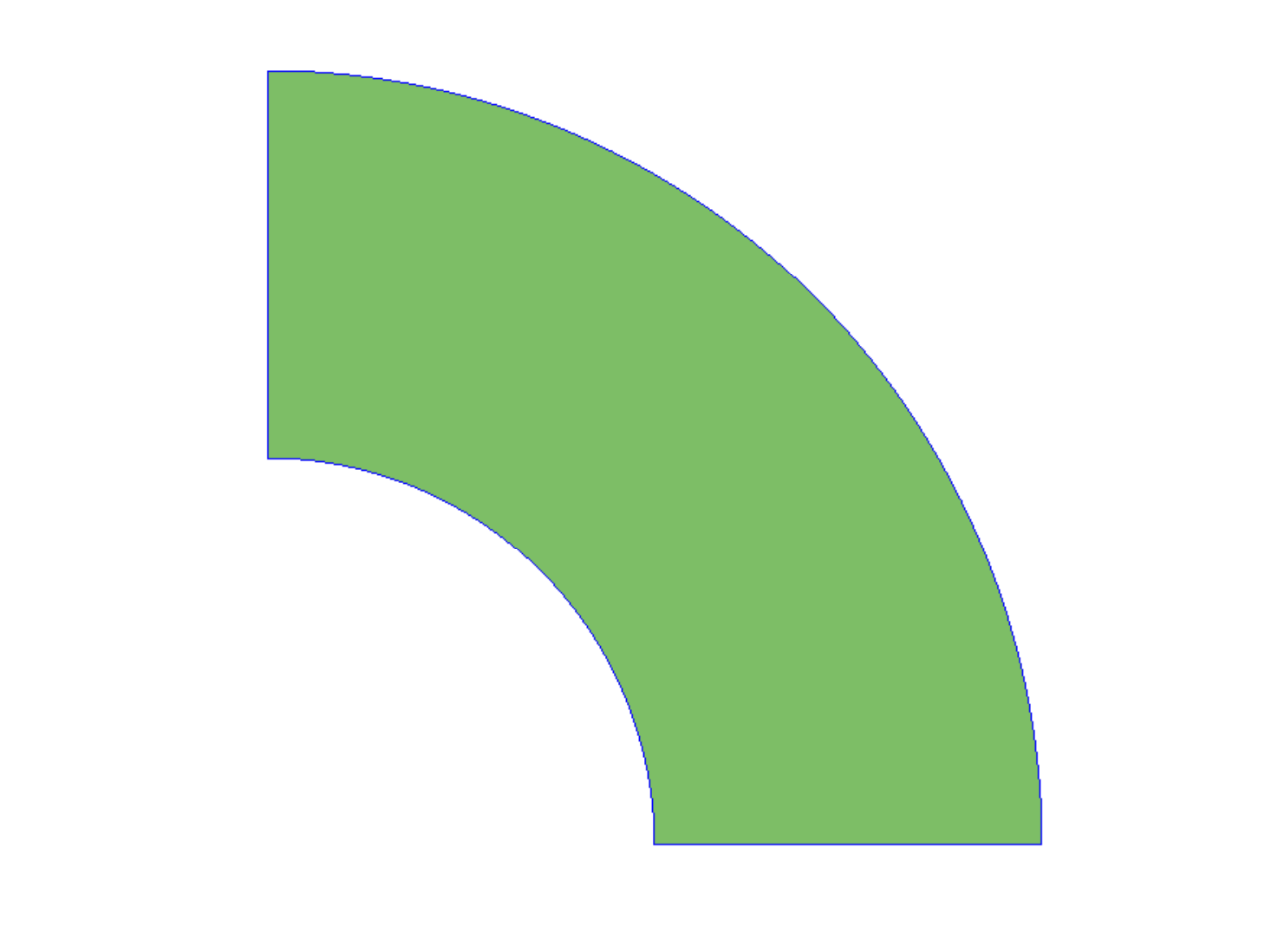}
\quad
\includegraphics[trim=0cm 1cm 2cm
1cm,angle=0,width=.45\columnwidth]{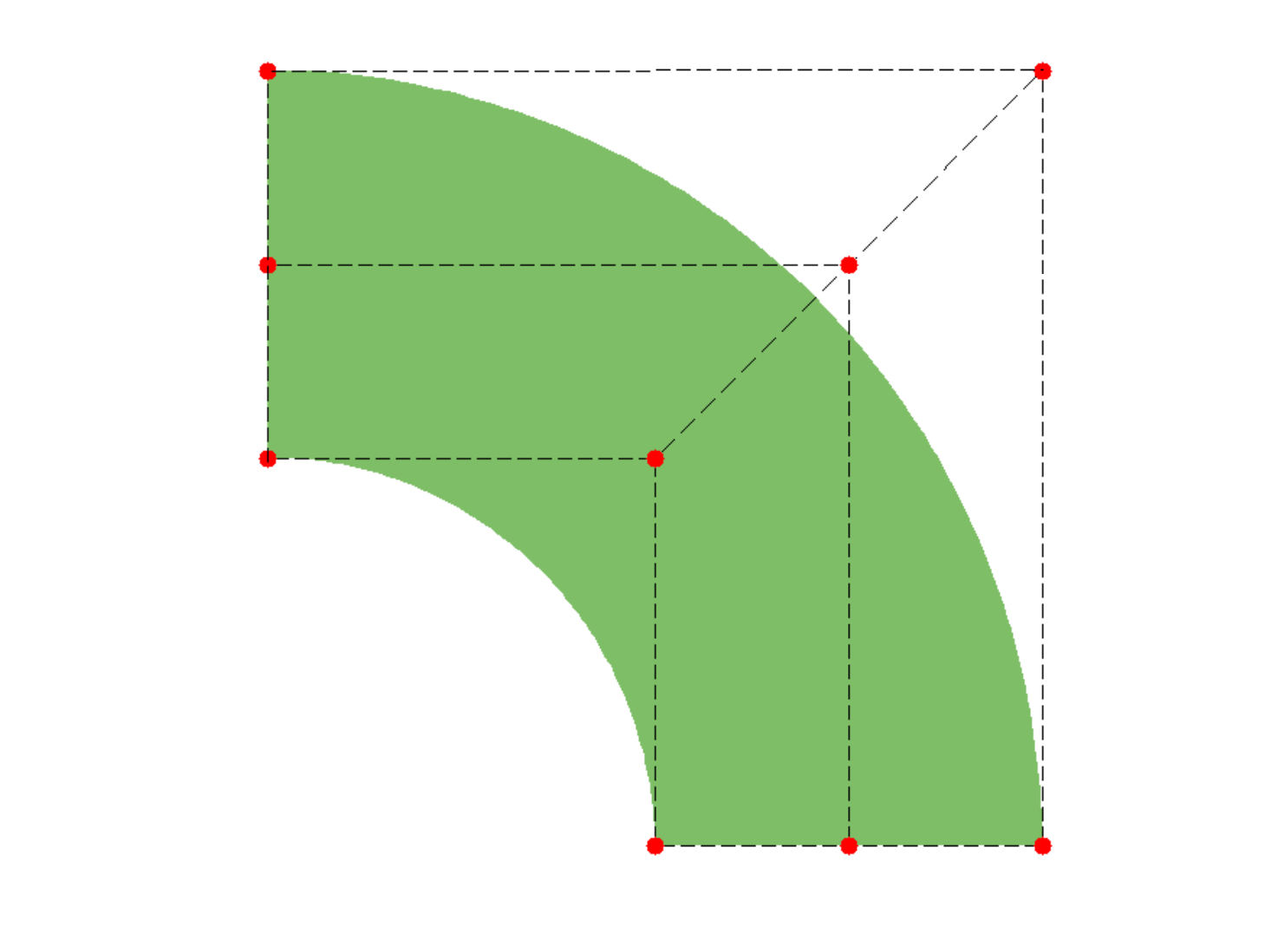}
\caption {Coarsest mesh (one single element) of the domain of problem~\ref{sec:problem3} and its control net.}
\label{fig:quartercylinder_coarse}
\end{figure}
%

\begin{figure}[!ht]
\centering
\includegraphics[angle=0,width=.32\columnwidth]{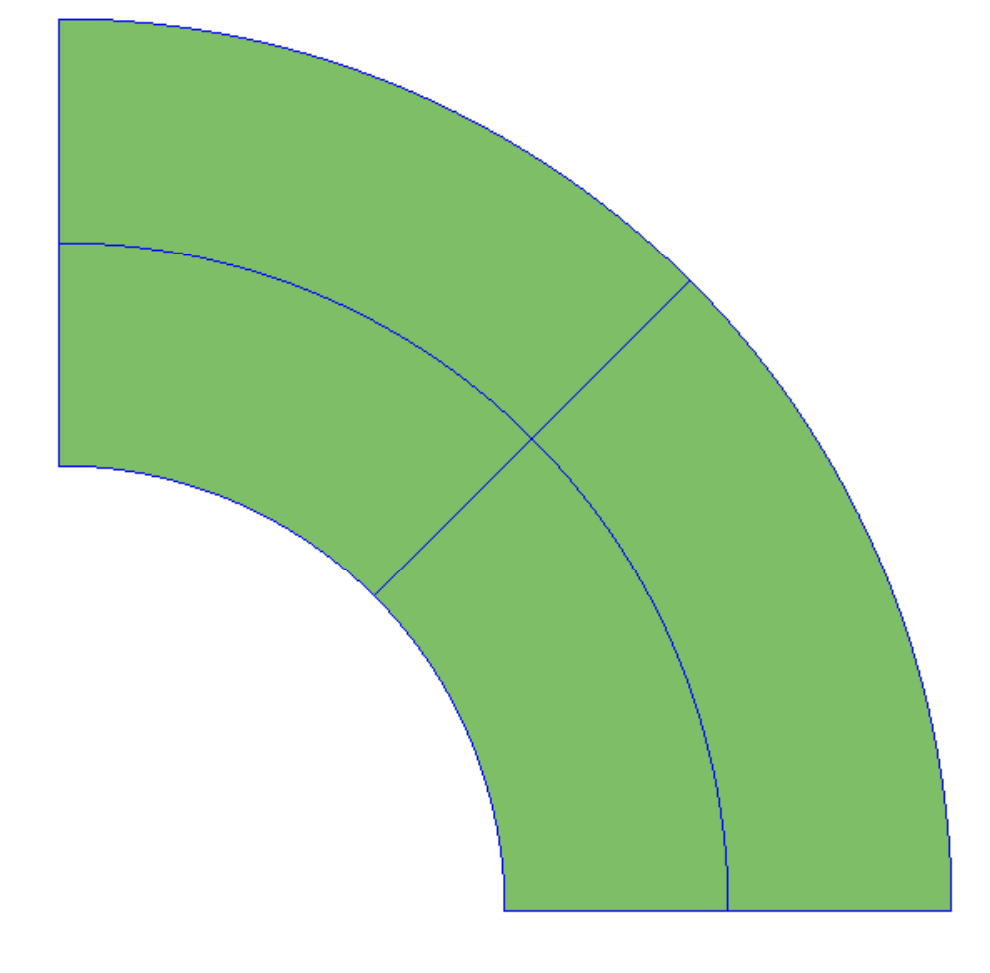}
\includegraphics[angle=0,width=.32\columnwidth]{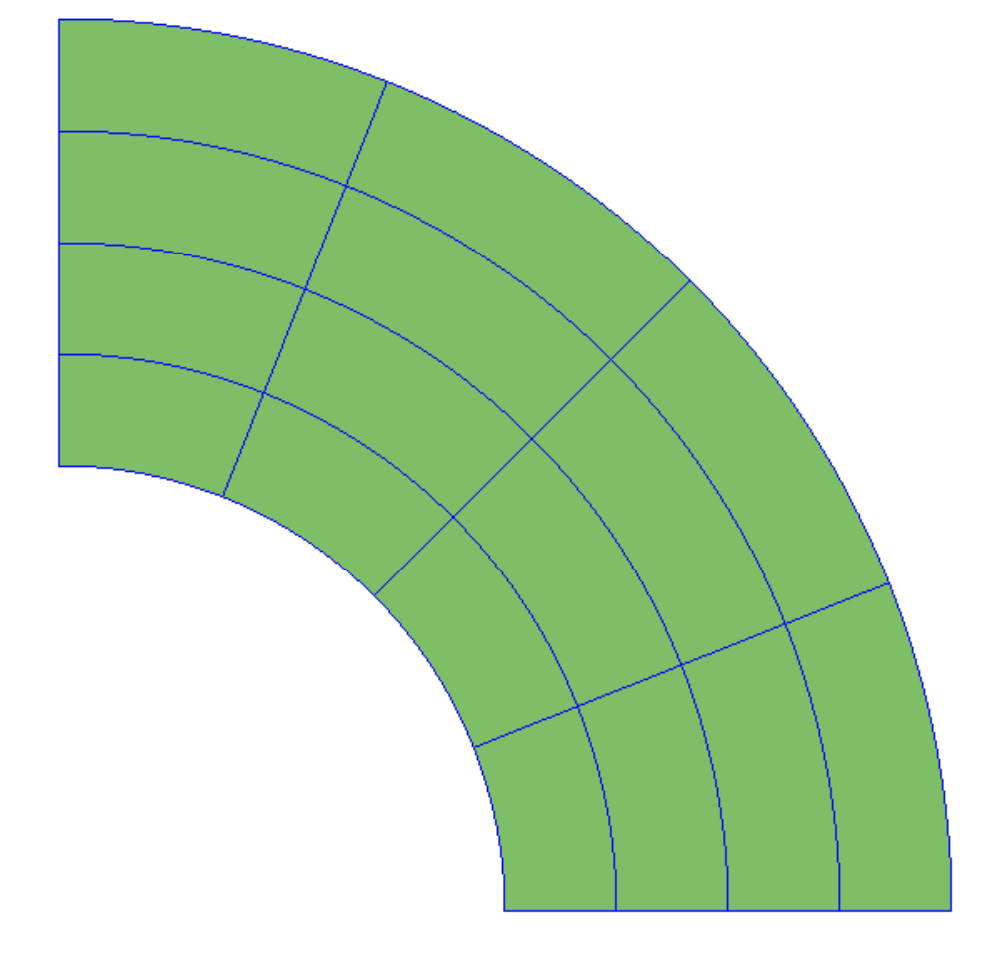}
\includegraphics[angle=0,width=.32\columnwidth]{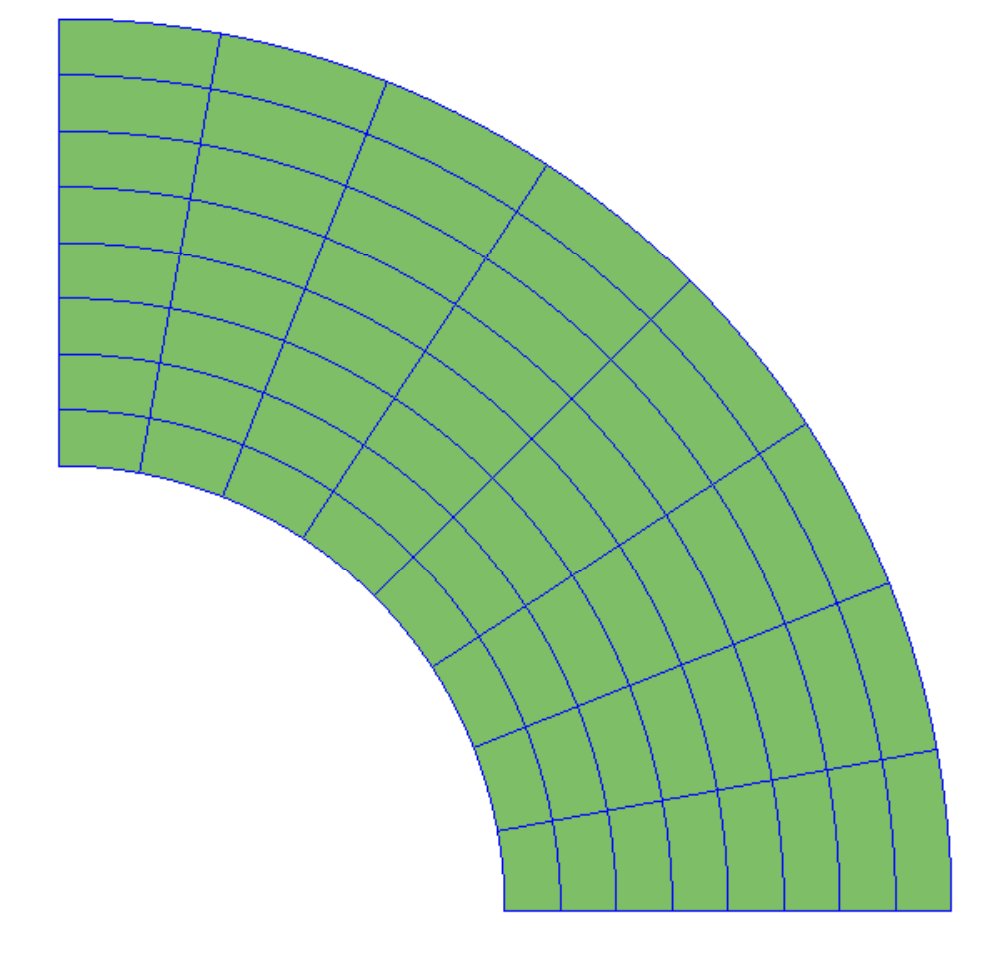}
\caption {Meshes ($2\times2$, $4\times4$, and $8\times8$) for the quarter circular annulus test cases.}
\label{fig:quartercylinder_meshes}
\end{figure}

\begin{table}[H]
\centering
\begin{tabularx}{\textwidth}{CCCCCCCCCC}
\toprule
$i$       & 1&2&3&4&5&6&7&8&9 \\
\hline
$\vm{P}_i$ &(0,1)&(1,1) &(1,0)&(0,1.5) &(1.5,1.5)&(1.5,0)& (0,2)&(2,2)&(2,0)\\
\hline
$w_i$&1&$1/\sqrt{2}$&1&1&$1/\sqrt{2}$&1&1&$1/\sqrt{2}$&1\\
\bottomrule
\end{tabularx}
\caption{Control points and weights for the coarsest mesh of the computational domain of
 problem \ref{sec:problem3}. The corresponding knot vectors are $\Xi  = \left\{ {0,0,0,1,1,1} \right\}, \mathcal{H}=\left\{0,0,0,1,1,1\right\}$.}
\label{tab:annulus_data}
\end{table}

\begin{table}[H]
\centering
\begin{tabularx}{\textwidth}{CCCC}
\toprule
&\multicolumn{3}{c}{Mesh}\\ 
\cline{2-4}
{Degree}      & $2\times2$  & $4\times4$ & $8\times8$ \\
\midrule
p=2      	 & 0.06920	 & 0.01682   & 0.00439\\
p=3       &	0.02575 &	0.00423 & 0.00052\\
p=4      	& 0.01567	  & 0.00107	 & 6.60E-05\\
p=5      	& 0.00713    &	0.00026 & 8.14E-06\\
\bottomrule
\end{tabularx}
\caption{$H^1$ error of the thermal quarter annulus problem,
using linear basis function in micro space and the $H^1$ micro refinement strategy.}
\label{tab:annulus_errH1_H1strategy}
\end{table}

\begin{figure}[!ht]
\begin{center}
\includegraphics[trim= 2cm 1cm 2cm 0cm,
angle=0,width=.45\columnwidth]{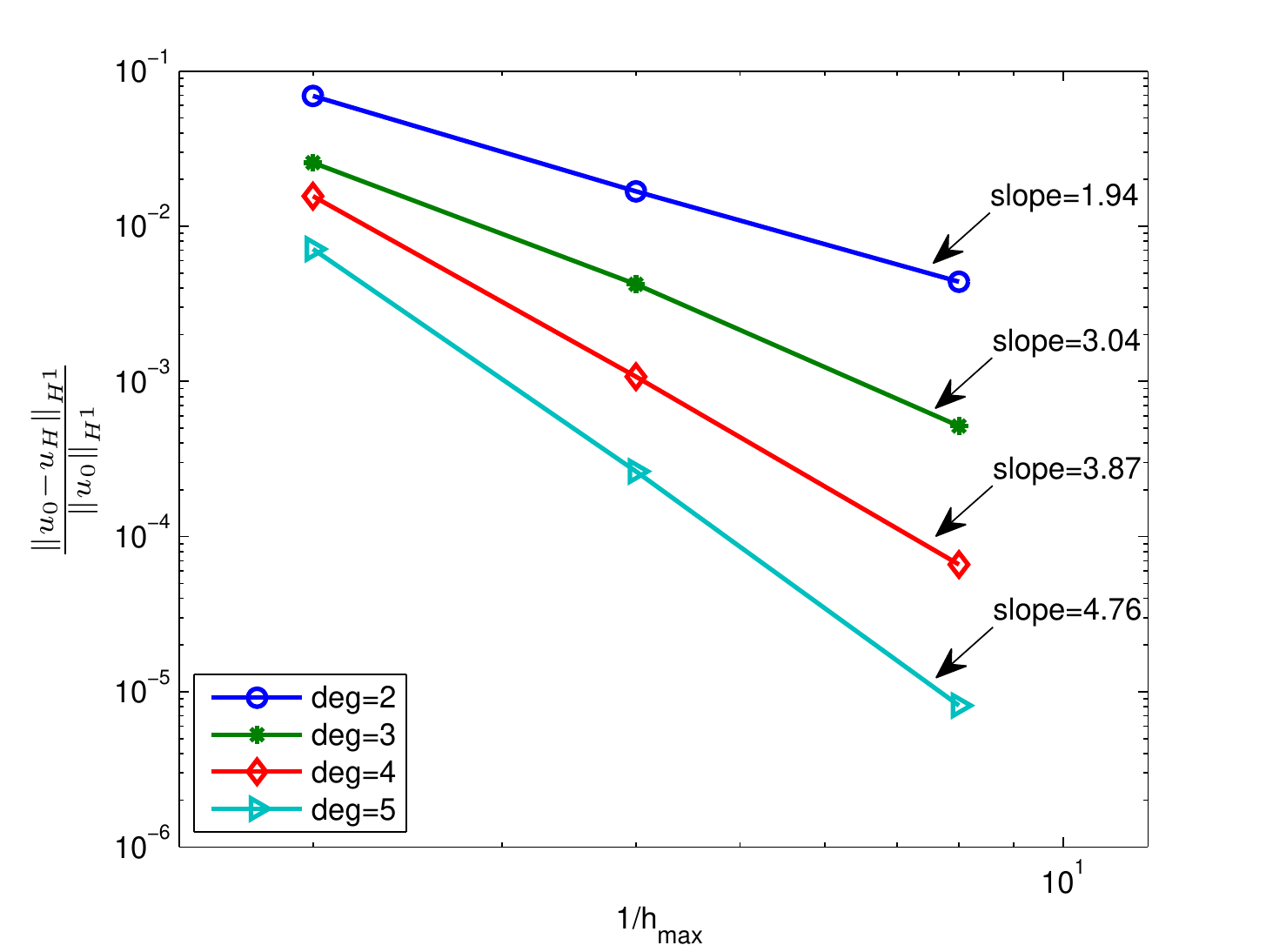}
\caption{$H^1$ error of the thermal quarter annulus problem,
using linear basis function in micro space and the $H^1$ micro refinement strategy.}
\label{fig:annulus_errH1_H1strategy}
\end{center}
\end{figure}


\begin{table}[h]
\centering
\begin{tabularx}{\textwidth}{CCCC}
\toprule
&\multicolumn{3}{c}{Mesh}\\ \cline{2-4}
{Degree}      & $2\times2$  & $4\times4$ & $8\times8$ \\
\midrule
p=2  &	0.05895	 & 0.01509	   &  0.00400\\
p=3  &	0.02501	 & 0.00401 &	0.00049\\
p=4  &	0.01514 & 0.00102 &	6.39E-05\\
p=5  &	0.00699 & 0.00026 &	7.91E-06\\
\bottomrule
\end{tabularx}
\caption{$L^2$ error of the thermal quarter annulus problem,
using linear basis function in micro space and the $H^1$ micro refinement strategy.}
\label{tab:annulus_errL2_H1strategy}
\end{table}

\begin{figure}[!ht]
\centering
\includegraphics[trim= 2cm 0cm 2cm 0cm,
angle=0,width=.45\columnwidth]{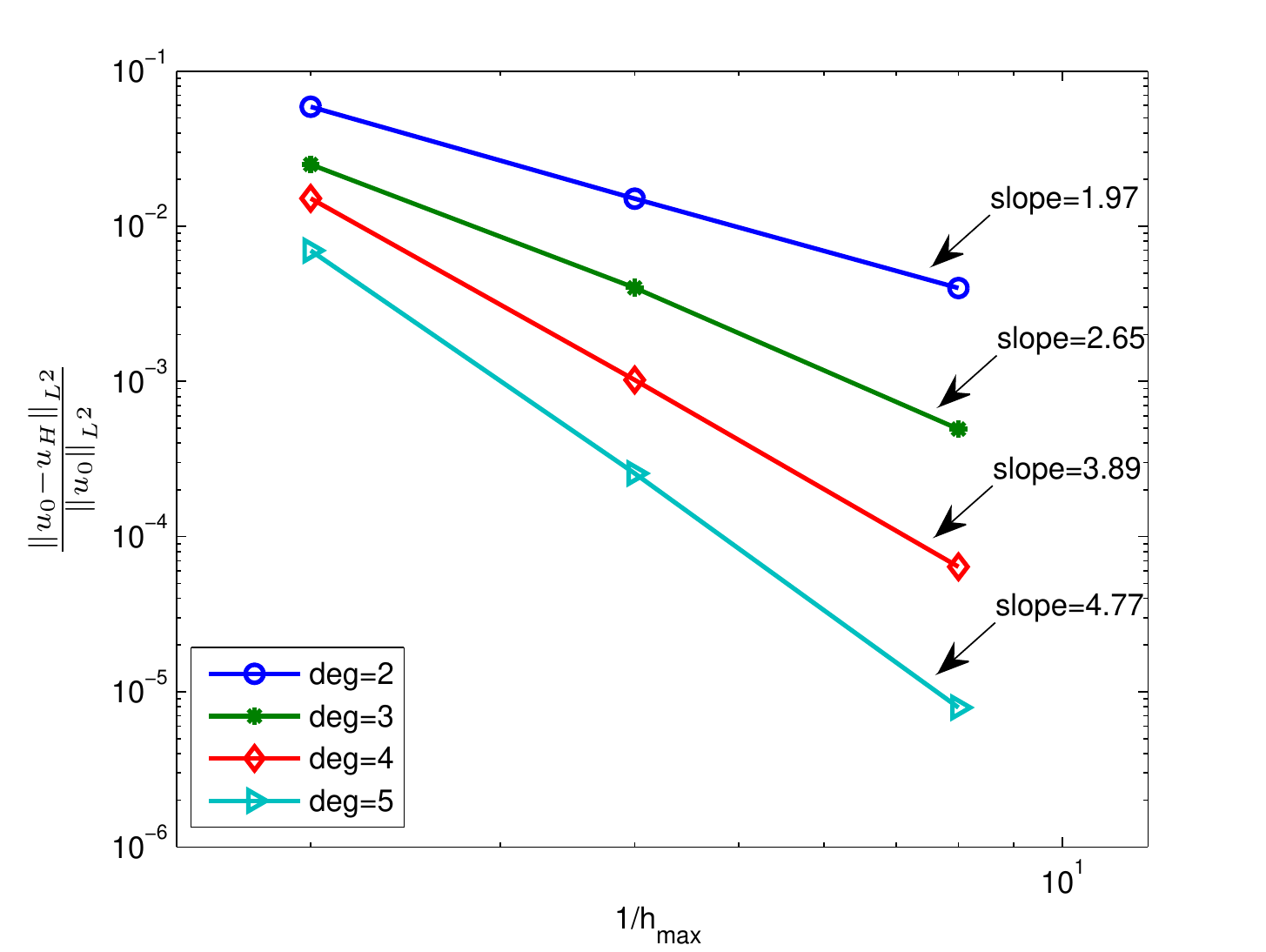}
\caption{$L^2$ error of the thermal quarter annulus problem,
using linear basis function in micro space and the  $H^1$ micro refinement strategy.}
\label{fig:annulus_errL2_H1strategy}
\end{figure}

We emphasize that in FE-HMM approach, it is very hard to use high order elements for problems with curved boundary due to the requirement
of the mesh quality: the partition of the domain must be fine
enough so that the curved boundary can be approximately good
enough. As a consequence, the computational cost will be
extremely high, if not  prohibited.

Next, we use high order NURBS basis functions in both micro and macro spaces. Macro space utilizes NURBS of degree ranging from $2$ to $5$, while in micro space quintic NURBS is used. Here, we follow the $L^2$ micro refinement strategy. The results in Tabs. \ref{tab:annulus_errH1_H1strategy_micdeg5}, \ref{tab:annulus_errL2_H1strategy_micdeg5}, and Figs. \ref{fig:annulus_errH1_H1strategy_micdeg5}, \ref{fig:annulus_errL2_H1strategy_micdeg5} show that this approach greatly outperforms the
standard FE-HMM in terms of  convergence rate, accuracy as well as efficiency.
In Tab. \ref{tab:annulus_time_H1strategy_micdeg5} and \Fref{fig:annulus_time_H1strategy_micdeg5}, the results also confirm the nearly optimal computational cost given by  Eq. \ref{eq:totalcost}.

It is worthy noting that for problems with curved boundaries as this example, it is very hard for FE-HMM to obtain accuracy of magnitudes $10^{-6}$, $10^{-7}$, $10^{-8}$ in $H^1$ error or $10^{-7}$, $10^{-8}$, $10^{-9}$ in $L^2$ error, because of the high number of macro elements required to approximate the boundary. Furthermore, elevating the degree of basis functions in standard FEM is, albeit possible, not a flexible task.  Here, in IGA-HMM with high order NURBS basis functions in both macro and micro space we can obtain these high accuracies very easily: about 1000 seconds with mesh $8 \times 8$  ($p=5$) for the accuracy of $10^{-6}$ in $H^1$ error  and $2 \times 10^{-7}$ in $L^2$ error; for accuracy of $2.7 \times 10^{-8}$ in $H^1$ error and $2.7 \times 10^{-9}$, it takes about 9000 seconds with a mesh of $16 \times 16$  ($p=5$). \\

\begin{table}[!ht]
\centering
\begin{tabularx}{\textwidth}{CCCCC}
\toprule
&\multicolumn{4}{c}{Mesh}\\ \cline{2-5}
{Degree} & $2\times2$ & $4\times4$ & $8\times8$ & $16\times16$ \\
\midrule
p=2&0.04501	&0.00724&	0.00164	&0.00040\\
p=3&0.00450&	0.00116	&1.13E-04&	1.33E-05\\
p=4&0.00241&	2.72E-04	&1.04E-05&	5.79E-07\\
p=5&0.00012	&3.80E-05&	1.11E-06&	2.71E-08\\
\bottomrule
\end{tabularx}
\caption{$H^1$ error of the thermal quarter annulus problem,
using NURBS of degree 5 in micro space}
\label{tab:annulus_errH1_H1strategy_micdeg5}
\end{table}

\begin{figure}[ht]
\begin{center}
\includegraphics[trim= 2cm 1cm 2cm 0cm,
angle=0,width=.45\columnwidth]{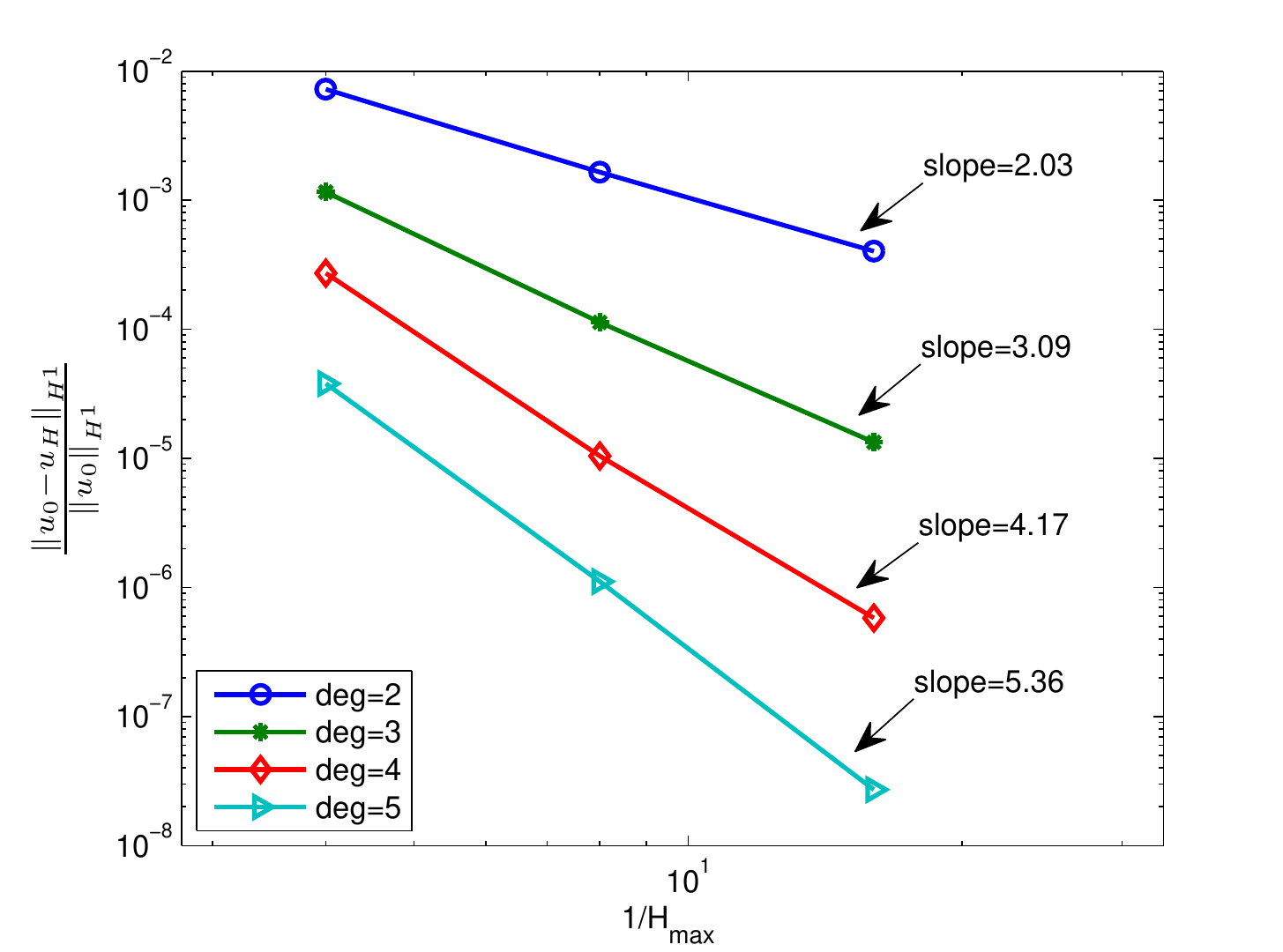}
\caption{$H^1$ error of the thermal quarter annulus problem, using NURBS of degree 5 in micro space.}
\label{fig:annulus_errH1_H1strategy_micdeg5}
\end{center}
\end{figure}

\begin{table}[ht]
\centering
\begin{tabularx}{\textwidth}{CCCCC}
\toprule
&\multicolumn{4}{c}{Mesh}\\ \cline{2-5}
{Degree} & $2\times2$ & $4\times4$ & $8\times8$ &$16\times16$ \\
\midrule
p=2  &	0.02784	 & 0.002084   &  0.00022 &2.58E-05\\
p=3  &	0.00211 &0.00041	&1.89E-05 &	1.08E-06 \\
p=4  &	0.00108    &1.10E-04         &1.84E-06	& 5.78E-08 \\
p=5  &	4.43E-05      &1.19E-05	     &2.08E-07&	2.71E-09\\
\bottomrule
\end{tabularx}
\caption{$L^2$ error of the thermal quarter annulus problem,
using NURBS of degree 5 in micro space.}
\label{tab:annulus_errL2_H1strategy_micdeg5}
\end{table}

\begin{figure}[ht]
\begin{center}
{\includegraphics[trim= 2cm 1cm 2cm 0cm,
angle=0,width=.45\columnwidth]{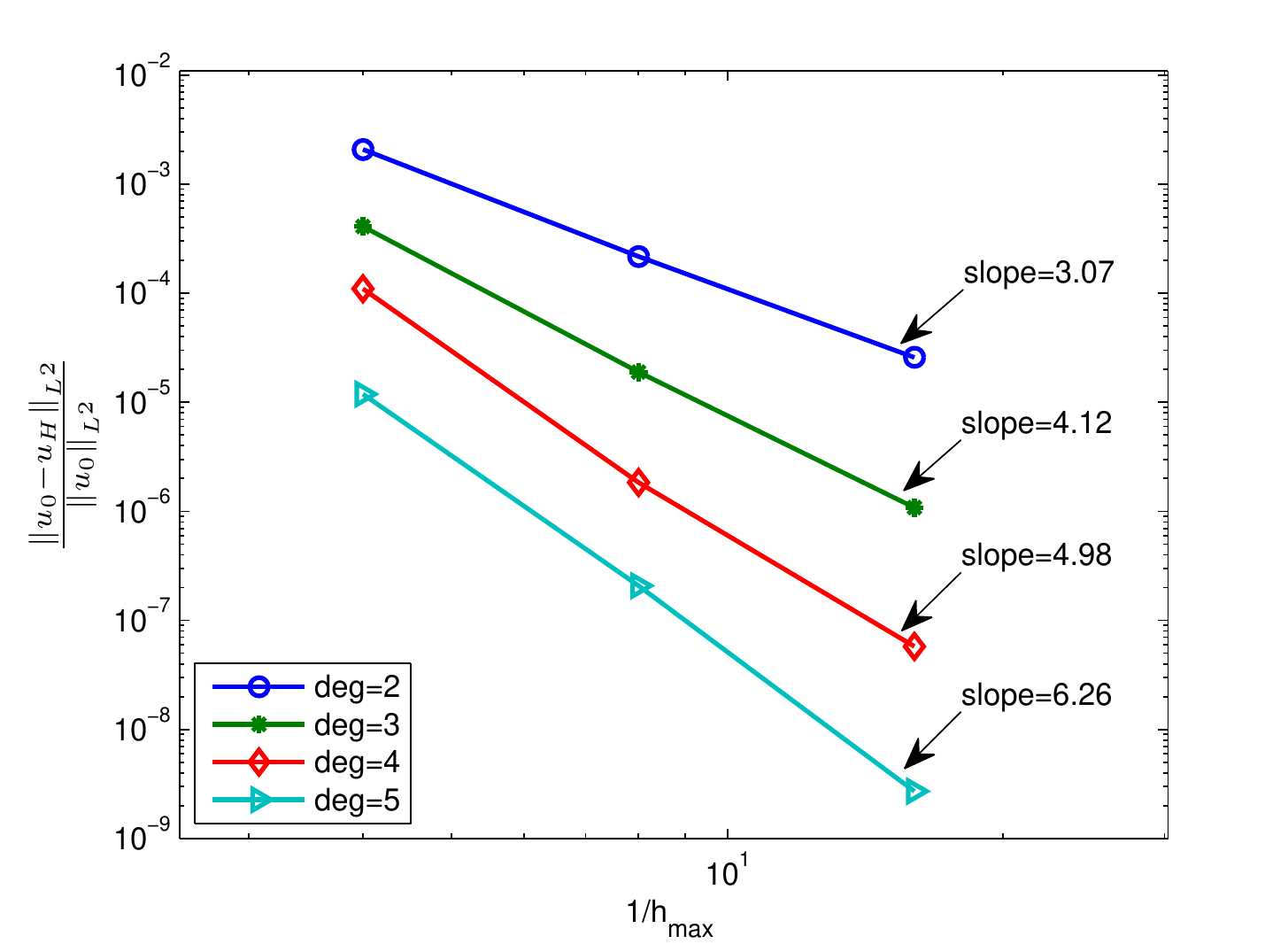}}
\caption{$L^2$ error of the thermal quarter annulus problem, using quintic NURBS in micro space.}
\label{fig:annulus_errL2_H1strategy_micdeg5}
\end{center}
\end{figure}

\begin{table}[H]
\centering
\begin{tabularx}{\textwidth}{CCCCC}
\toprule
&\multicolumn{4}{c}{Mesh}\\ \cline{2-5}
Degree & $2\times2$ & $4\times4$ & $8\times8$&$16\times16$ \\
\midrule
p=2&4.27    &	17.01 &  68.02     &	512.62\\
p=3&7.96    &	32.33 &	235.54 &	1638.44\\
p=4&12.48 &	 49.84   &	  374.17	& 2592.46\\
p=5&18.84	& 141.07 &	  956.76  &	9055.52\\
\bottomrule
\end{tabularx}
\caption{CPU time (s) solving the thermal quarter annulus problem,
using quintic NURBS in micro space.}
\label{tab:annulus_time_H1strategy_micdeg5}
\end{table}

\begin{figure}[H]
\begin{center}
{\includegraphics[trim= 2cm 1cm 2cm 1cm,
angle=0,width=.45\columnwidth]{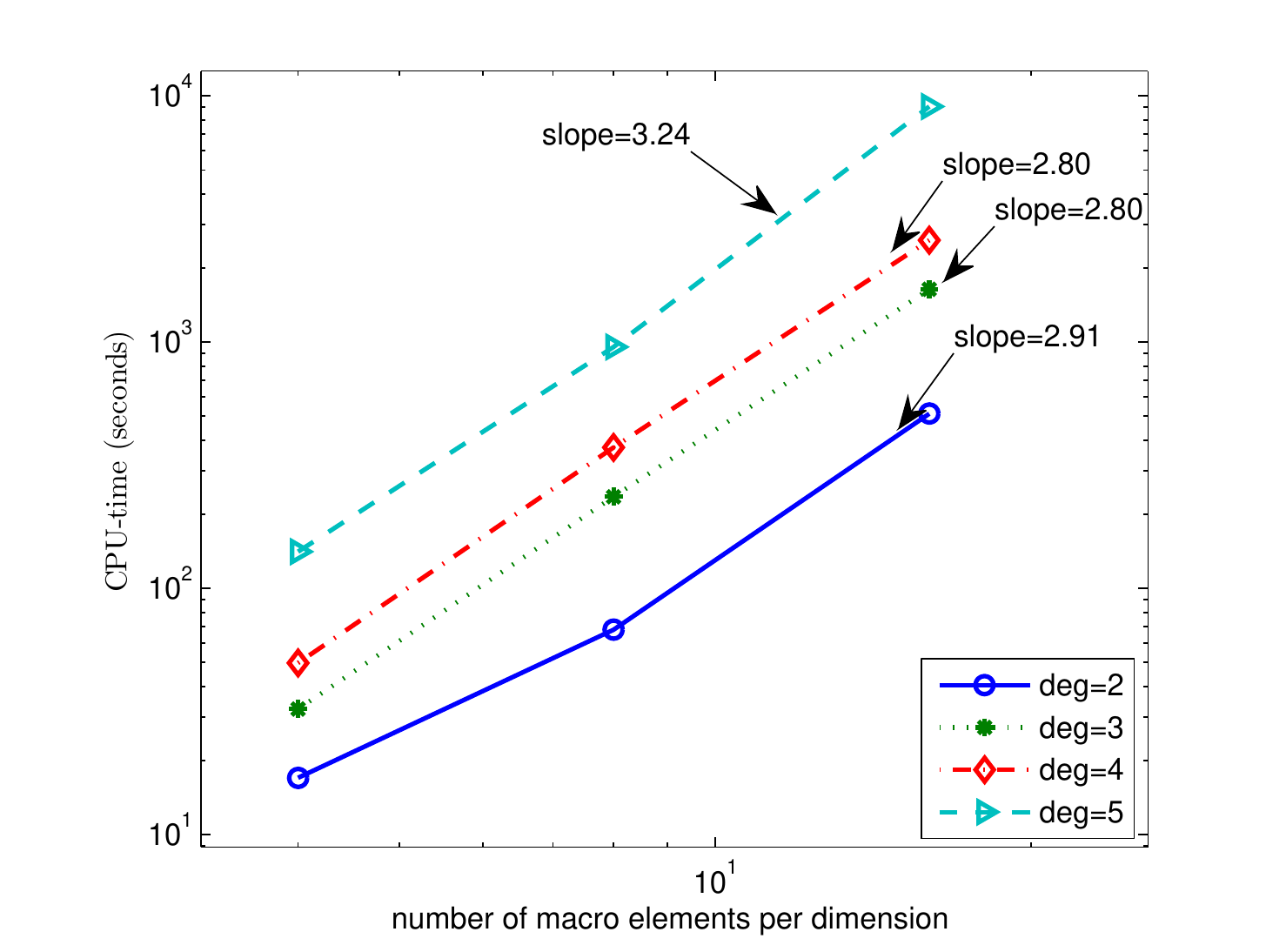}}
\caption{CPU time (s) solving thermal quarter annulus problem, using quintic NURBS in micro space.}
\label{fig:annulus_time_H1strategy_micdeg5}
\end{center}
\end{figure}

Hitherto, $C^0$ NURBS are used in the micro spaces. 
In what follows, we test the performance of the IGA-HMM when utilizing high order NURBS with $C^{p-1}$ continuity in  the macro space and
with $C^{q-1}$ continuity in the micro space. 
The errors in $H^1$ norm and $L^2$ norm are given in Tabs. \ref{tab:annulus_errH1_H1strategy_micdeg5_Cp} and \ref{tab:annulus_errL2_H1strategy_micdeg5_Cp}, respectively.
In this case, the performance of the method is worst than the performance when $C^0$ NURBS is used in micro spaces, and the convergence rate is not optimal, see Fig. \ref{fig:problem43_Cp_err_L2strategy}.
The reason for this phenomenon, to our knowledge, is that in the micro space the material is strongly heterogeneous as reflected in 
$a\left( {\frac{\mathbf{x}}{\varepsilon}} \right)$. As a consequence, high order $C^{q-1}$ continuity NURBS basis functions does not
capture well these heterogeneities as the $C^0$ ones, although the accuracy is still acceptable in comparison with that obtained with
the standard FE-HMM. 

\begin{table}[!ht]
\centering
\begin{tabularx}{\textwidth}{CCCCCC}
\toprule
&\multicolumn{4}{c}{Mesh}\\ \cline{2-6}
{Degree} & 4x4 & 8x8&16x16 & 32x32 & 64x64\\
\midrule
p=2  &	0.045027	 & 0.007264   &  0.001738 &4.67E-04&2.70E-4\\
p=3  &	0.004546 &0.001299	&2.7E-4 &	1.0E-04 & 1.0E-4 \\
\bottomrule
\end{tabularx}
\caption{$H^1$ error of the thermal quarter annulus problem, using quintic NURBS in micro space with $C^0$ continuity}
\label{tab:annulus_errH1_H1strategy_micdeg5_Cp}
\end{table}

\begin{table}[!ht]
\centering
\begin{tabularx}{\textwidth}{CCCCCC}
\toprule
&\multicolumn{4}{c}{Mesh}\\ \cline{2-6}
{Degree} & 4x4 & 8x8&16x16 & 32x32 & 64x64\\
\midrule
p=2  &	0.027884	 & 0.002179   &  0.000553 & 0.000221 & 0.000230 \\
p=3  &	0.002213 &0.00070 	& 0.000231 &	9.27E-05 & 9.42E-5 \\
\bottomrule
\end{tabularx}
\caption{$L^2$ error of the thermal quarter annulus problem, using quintic NURBS  in micro space with $C^0$ continuity}
\label{tab:annulus_errL2_H1strategy_micdeg5_Cp}
\end{table}

\begin{figure}[H]
\begin{center}
\subfigure{\includegraphics[width=.49\columnwidth]{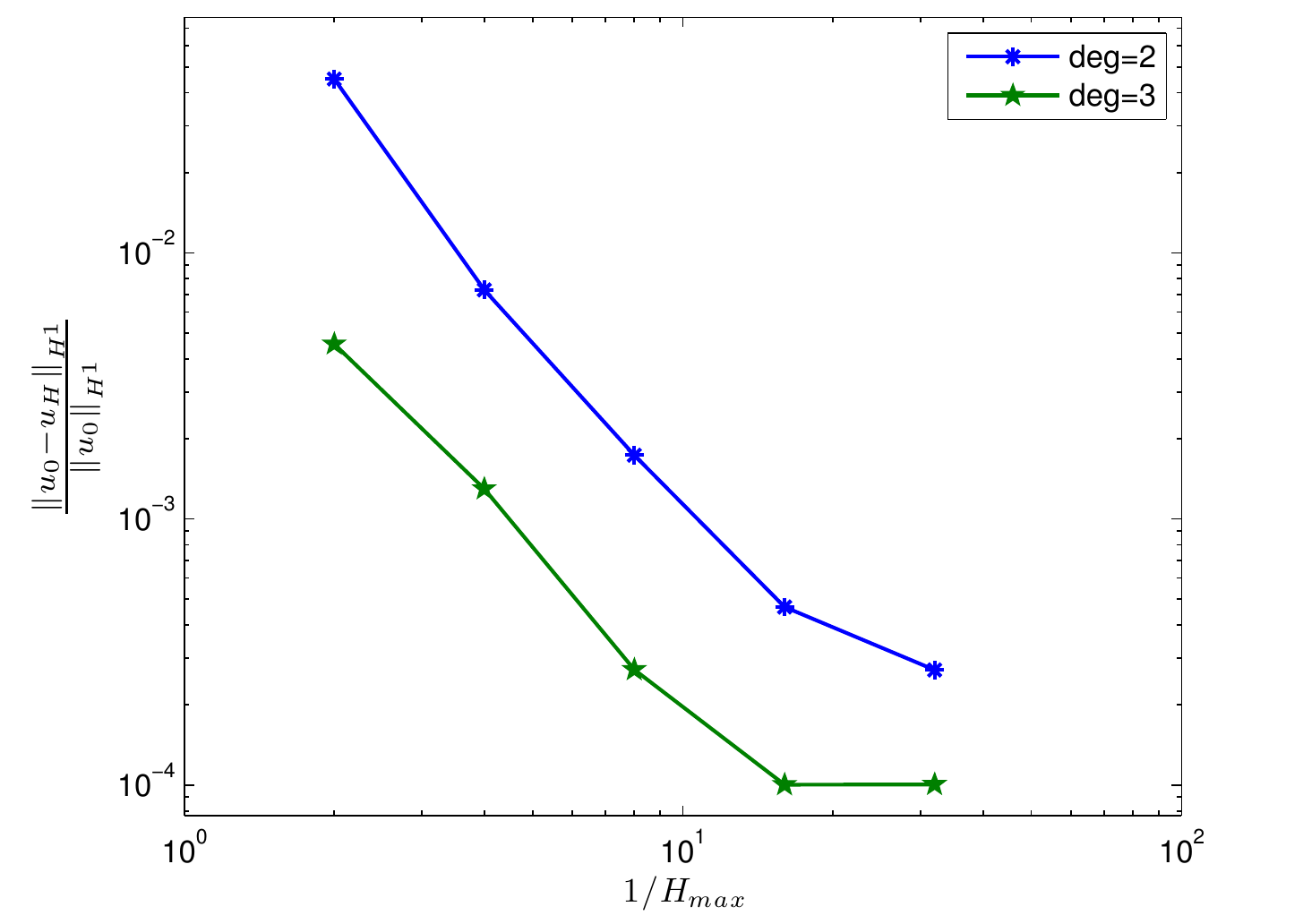}}
\subfigure{\includegraphics[width=.49\columnwidth]{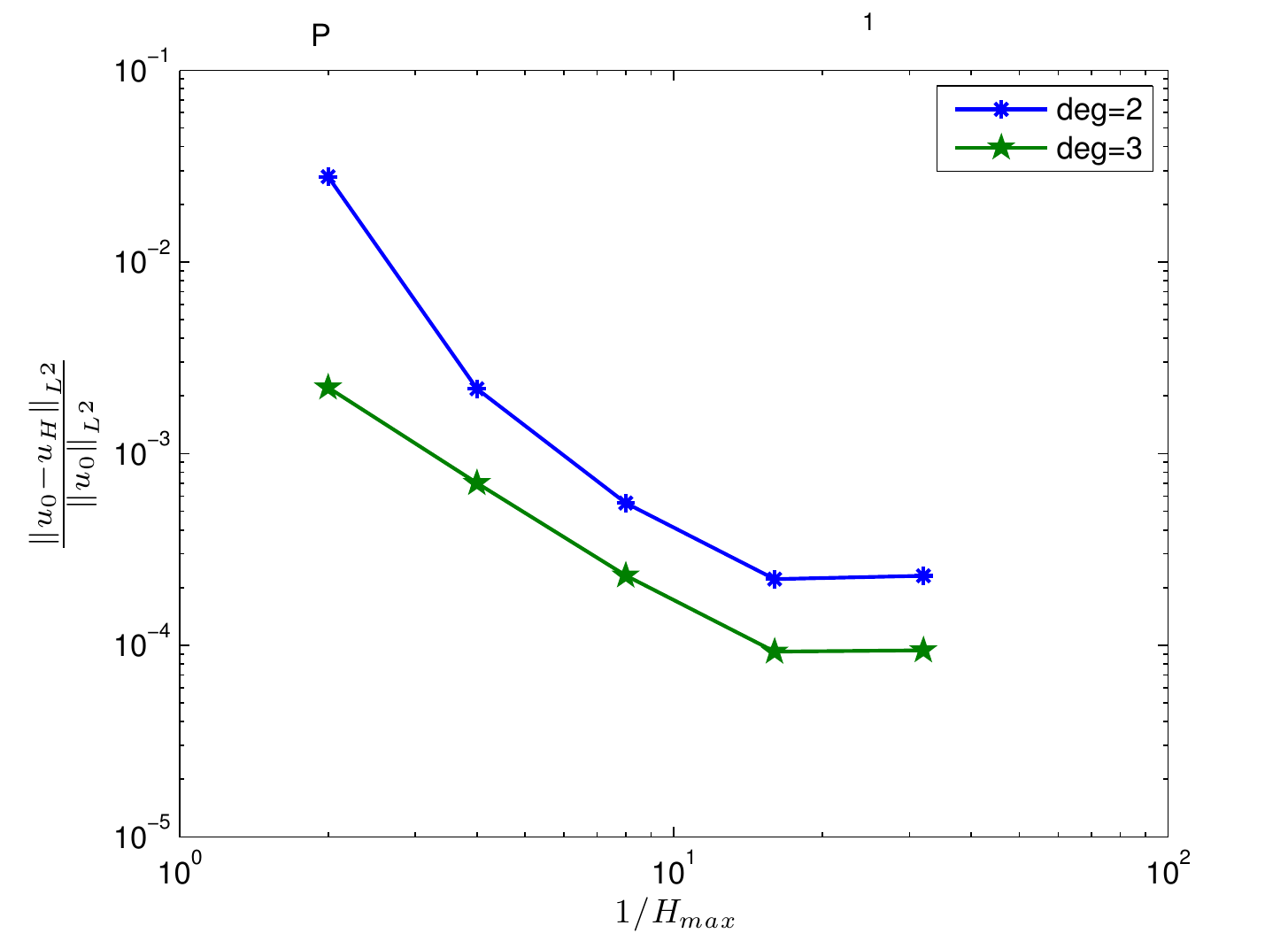}}
\end{center}
\caption{{\color{black}{$H^1$ and $L^2$ errors of the thermal quarter annulus problem, using quintic NURBS in the micro space with $C^{p-1}$ continuity.}}}
\label{fig:problem43_Cp_err_L2strategy}
\end{figure}

\section{Conclusions}
We have, for the first time, presented an efficient isogeometric analysis heterogeneous multiscale method (IGA-HMM) for elliptic homogenization problems. The method is capable of capturing the exact geometric representation and is very flexible regarding to refinement and degree elevation by using the NURBS basis functions in both macro and micro levels. As a result, the high-order IGA-HMM macroscopic and microscopic solvers are designed simply and effectively. A priori error estimates of the discretization errors and optimal micro refinement strategies were derived. The provided numerical results showed that the IGA-HMM achieves high order accuracy with  optimally convergence rate for both $L^2$ and $H^1$ micro refinement strategies as in the theory of the standard FE-HMM. 
We have also observed that $C^0$, not $C^{p-1}$, high order NURBS should be used at the microscale although this point deserves a further study.
The importance is that this approach is very similar to the standard FE-HMM, this means that the implementation is very simple, the coding framework of FE-HMM in \citep{AbdulleShort2009} can be re-used, while the accuracy and the flexibility are remarkably improved. These advantages make the IGA-HMM  an alternative method to solve  homogenization problems, beside some new methods have recently been  published \cite{li_efficient_2012,abdulle_reduced_2012}. The method shows potential in solving high order homogenization problems due to the arbitrary smoothness of the NURBS.

In order to apply the IGA-HMM for problems with a more higher complexity, a combination of IGA-HMM with techniques introduced in \cite{patera_reduced_2007} can be considered. In this way, not only the number of macro elements will be reduced but also that of micro problems. As a result, the computational cost will dramatically be minimized. This will be our forthcoming development.\\

\noindent\textbf{Acknowledgement} The support of the Vietnam National Foundation for Science and Technology Development (NAFOSTED) is gratefully acknowledged.


%
\section*{References}
\bibliography{References}
\bibliographystyle{plain}
\end{document}